\documentclass{notices}

\usepackage{amsmath,amssymb,amsthm}
\usepackage{graphicx}
\usepackage{tikz}
\usepackage{enumerate}
\usepackage{algorithm,algpseudocode}

\algdef{SE}[SUBALG]{Indent}{EndIndent}{}{\algorithmicend\ }%
\algtext*{Indent}
\algtext*{EndIndent}

\usetikzlibrary{arrows,decorations.markings}
\tikzstyle{vecArrow} = [thick, decoration={markings,mark=at position
   1 with {\arrow[semithick]{open triangle 60}}},
   double distance=1.4pt, shorten >= 5.5pt,
   preaction = {decorate},
   postaction = {draw,line width=1.4pt, white,shorten >= 4.5pt}]

\newcommand{\eps}{\epsilon}

\newcommand{\rmd}{\mathrm{d}}

\newcommand{\rmr}{\mathrm{r}}

\newcommand{\Nb}{\mathbb{N}}
\newcommand{\Rb}{\mathbb{R}}
\newcommand{\Sb}{\mathbb{S}}
\newcommand{\Vb}{\mathbb{V}}

\newcommand{\Ac}{\mathcal{A}}
\newcommand{\Bc}{\mathcal{B}}

\newcommand{\Hc}{\mathcal{H}}
\newcommand{\Ic}{\mathcal{I}}
\newcommand{\Kc}{\mathcal{K}}
\newcommand{\Lc}{\mathcal{L}}
\newcommand{\Nc}{\mathcal{N}}
\newcommand{\Pc}{\mathcal{P}}
\newcommand{\Sc}{\mathcal{S}}

\newcommand{\Xc}{\mathcal{X}}
\newcommand{\Yc}{\mathcal{Y}}

\newcommand{\As}{\mathsf{A}}
\newcommand{\Bs}{\mathsf{B}}
\newcommand{\Gs}{\mathsf{G}}
\newcommand{\Ps}{\mathsf{P}}
\newcommand{\Qs}{\mathsf{Q}}
\newcommand{\Rs}{\mathsf{R}}
\newcommand{\Ss}{\mathsf{S}}
\newcommand{\Us}{\mathsf{U}}
\newcommand{\Vs}{\mathsf{V}}
\newcommand{\Ys}{\mathsf{Y}}

\newcommand{\ran}{\mathrm{ran}}
\newcommand{\diag}{\mathrm{diag}}

\theoremstyle{definition}\newtheorem{theorem}{Theorem}
\theoremstyle{definition}
\theoremstyle{definition}\newtheorem{definition}{Definition}
\theoremstyle{definition}\newtheorem{proposition}{Proposition}

\graphicspath{{./figures/}}

\title{
Low-rank approximation for multiscale PDEs
}

\author{
  Ke Chen
  \affil{
    Ke Chen is a postdoctoral researcher in the Department of Mathematics at University of Texas at Austin. His email address is ke.chen@austin.utexas.edu.
    }
  \and
  Shi Chen
  \affil{
    Shi Chen is a graduate student in the Department of Mathematics at University of Wisconsin, Madison. His email address is schen636@wisc.edu.
   }
  \and
  Qin Li
  \affil{
    Qin Li is an associate professor in the Department of Mathematics and the Wisconsin Institutes for Discovery at University of Wisconsin, Madison. Her email address is qinli@math.wisc.edu.
   }
  \and
  Jianfeng Lu
  \affil{
    Jianfeng Lu is a professor in the Department of Mathematics, the Department of Physics and the Department of Chemistry at Duke University. His email address is jianfeng@math.duke.edu.
   }
  \and
  Stephen J. Wright
  \affil{
    Stephen J. Wright is a professor in the Department of Computer Sciences at University of Wisconsin, Madison. His email address is swright@cs.wisc.edu.
   }
}

\begin{document}

\maketitle

\section{Introduction}
Multiscale phenomena are ubiquitous, with applications in many physical sciences and engineering fields:  aerospace, material sciences, geological structure analysis, and many others.
The different scales often have different physics, which entangle to produce complicated nonlinearities.
Partial differential equations (PDEs) are often used to model these problems, with different scales captured in the coefficients and functions that define the PDE.
These PDE models are challenging to compute directly, so analysis and algorithms specifically targeted to multiscale problems have been developed and investigated.
Following convention, we focus in this review on problems with two distinct scales, with a small positive parameter $\epsilon$ capturing the ratio between the small and large scale.

Though modern multiscale analysis dates back to asymptotic PDE analysis that was seen already in Hilbert and Poincar\'e expansions early last century (see review in~\cite{MR2382139}), the impetus for computations involving multiscale PDEs  came largely from the US Department of Energy (DOE) National Labs within the ASCI (Advanced Strategic Computing Initiative)~\cite{Ho:2009multiscale}.
Since that time, analysis and computation in multiscale PDEs have taken different paths.
Analysis has tended to follow a single ``universal'' strategy, passed down from tradition.
The equation is decomposed into several levels according to asymptotic expansions involving the scale parameters, with the subequation at each level representing physics at a single scale, and subequations at the finer level feeding information to those at the coarser level.
This analytical machinery has been used to treat multiscale PDEs arising from such varied backgrounds as kinetic theory, semi-classical quantum systems, and homogenization of composite materials, among others \cites{MR2382139, MR2830582}.

On the computational side, strategies for handling multiscale PDEs are more varied.
Problems are usually handled by specifically designed solvers.
One class of solvers called  {\em asymptotic-preserving} schemes~\cite{MR3645390} are designed to preserve asymptotic limits of kinetic equations.
These schemes usually contain some component of macro-solvers and micro-solvers, integrated in a clever way to reveal different structures in different regimes.
Another class of solvers called {\em numerical homogenization} methods~\cites{MR2830582,MR3971243,MR1979846} usually target elliptic and parabolic equations in which the coefficients that represent media have oscillatory elements.
These methods usually consist of offline and online stages, with either the homogenized media or the representative basis functions being prepared in the offline stage.

Why are most numerical schemes for multiscale PDEs equation-specific despite the analytical tools being largely unified?
This intriguing question has motivated our investigations into devising a universal numerical strategy for solving multiscale PDEs.
While the approach is yet to be developed fully, we believe that our progress on this issue is of wide interest, and this article surveys our progress to date.
Crucially, our approach exploits the low-rank structure present in discretizations of multiscale PDEs.

To demonstrate the fundamental idea, we consider the following problem:
\begin{equation}\label{eqn:general}
\Lc^\eps u^\eps = f\,,
\end{equation}
where $\Lc^\eps$ is a linear partial differential operator that depends explicitly on the small parameter $\eps$, while $f$ represents the boundary condition or the source term, which is assumed to have no dependence on $\eps$.
Multiscale problems that can be formulated in this way include elliptic equations with highly oscillating media and the neutron transport equation with small Knudsen number.
Due to the $\eps$-dependence of the operator, the solution $u^\eps$ inherits structures at both fine and coarse scales.

An asymptotic limit is revealed by multiscale analysis using asymptotic expansions as $\eps \to 0$.
In this limit, the oscillation at the fine scale is fast and the detailed oscillation pattern no longer matters --- only macroscopic quantities are relevant.
Formally, writing the homogenization limit as
\begin{equation}\label{eqn:general_homo}
\Lc^\ast u^\ast = f \,,
\end{equation}
we have
\begin{equation}
\| u^\eps-u^\ast \| \to 0 \textrm{ as }\eps\to0 \,.
\end{equation}
The norm of the approximation error depends heavily on the particular equation at hand.

The numerical challenge in solving~\eqref{eqn:general} is that many degrees of freedom may be needed.
Na\"ive finite element or finite difference methods would require mesh size $h\ll\eps$ to resolve fine-scale structure of the solution at the $\eps$ level.
For a problem on $\mathbb{R}^d$, the discretized system would therefore have $O(\eps^{-d})$ degrees of freedom, leading to prohibitive computational and memory cost for small $\eps$.
From an application perspective, it often suffices to characterize the solutions on the macroscopic level, where oscillations at the $\eps$ scale are largely absent.
This property raises the question of whether we can obtain an approximate solution of this type using only $O(1)$ degrees of freedom.
If we know how to derive \eqref{eqn:general_homo}, we can simply solve for $u^\ast$, which has the required macroscopic properties, and typically requires a discretization with $O(1)$ degrees of freedom.
Often, though, the limiting equation \eqref{eqn:general_homo} and its solution $u^\ast$ are difficult to find explicitly, even when it is possible to establish their existence.
These difficulties have led researchers to propose problem-specific solutions.

\begin{figure}[tb]
	\centering
 	\includegraphics[width = 0.5\textwidth]{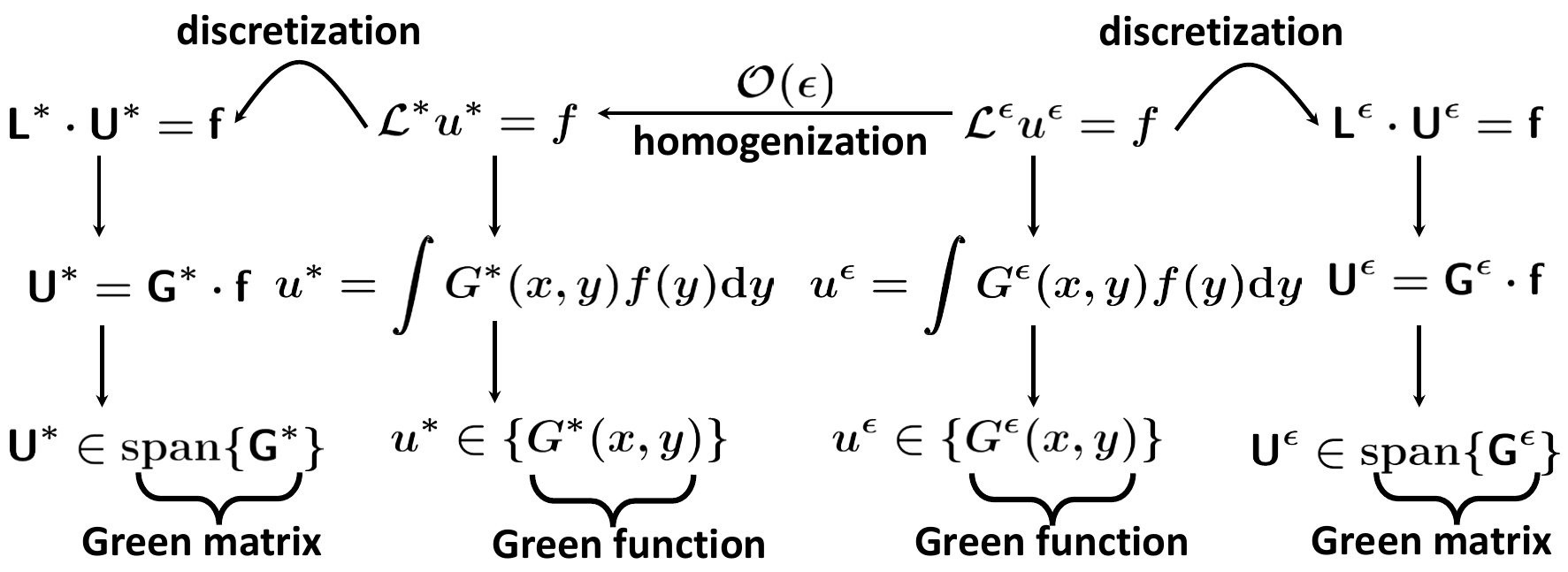}
	\caption{PDEs with small parameters have homogenized limits, meaning the solutions to the original PDEs can be well-approximated by the solutions to the limiting equations. While analytically the two solution spaces are ``close", the original equation requires many more degrees of freedom to solve numerically than its limiting counterpart. The numerical Green's matrix is intrinsically low rank.}
	\label{fig:low_rank}
\end{figure}

We believe that a universal approach for finding the large-scale solution can be devised, and that exploitation of the low-rank structure of the solution space is the key to developing such an approach.
As suggested above, the Green's matrix $\mathsf{G}^\eps$ (the discretized Green's function on fine grids) for the multiscale system~\eqref{eqn:general}  requires dimension  $O(\epsilon^{-d})$ to represent the underlying Green's function accurately.
However, if a limiting system such as~\eqref{eqn:general_homo} exists, this limiting system can be well-represented numerically by $\mathsf{G}^\ast$, a Green's matrix with dimension only $O(1)$. This phenomenon suggests the system can largely be ``compressed'' and hence is of low rank; see illustration in Figure~\ref{fig:low_rank}. In the language of  numerical linear algebra, this transition amounts to performing a truncated singular value decomposition  (SVD) of $\mathsf{G}^\eps$ to obtain $\mathsf{G}^\ast$.

If we obtain the truncated SVD of the matrix $\mathsf{G}^\eps$ by starting with a {\em full} SVD, the resulting algorithms would be impractical because of the large dimension of the matrix and the expense of preparing and storing the full matrix $\mathsf{G}^\eps$ and computing its SVD.
Several new linear algebra solvers take a quite different approach. Instead of accessing the full matrix, these new solvers merely require computation of matrix-vector products, involving the target matrix and several randomly selected vectors (typically vectors with Gaussian i.i.d. entries).
Translated to the PDE solver setting, these matrix-vector multiplications amount to computing numerical solutions to PDEs with some random source terms, a task that may be practical if the number of such operations required is modest.
The randomized SVD (rSVD) approach is one method of this type.
It is equipped with a thorough analysis and achieves optimality in terms of computational efficiency.
We make use of this method in the techniques described in the remainder of this article.

The main theme of our article, then, is the use of randomized SVD solvers to exploit the low-rank features of multiscale PDEs.
We will describe two strategies both of which are divided into ``offline" and ``online" stages.
The offline stage sees the preparation of either the solution space or the boundary-to-boundary map used in the domain decomposition, while the online stage singles out the specific solution for the given source $f$.
The two strategies are described in Section~\ref{sec:learning_basis} and~\ref{sec:Schwarz}, respectively.
In Section~\ref{sec:manifold}, we present the nonlinear extension utilizing manifold learning algorithms for reconstructing the low-rank features of the solution manifold.
Prior to these discussions, we describe in Section~\ref{sec:models} two algorithm classes --- asymptotic preserving and numerical homogenization --- for identifying the asymptotic limits of multiscale problems.
As examples, we use the multiscale radiative transfer equation (RTE) and the elliptic equation with rough media.
Section~\ref{sec:framework} explores the two main elements of our approaches: the numerical low-rank feature of multiscale PDEs and the randomized SVD solver for efficient reconstruction of low-rank operator/spaces.
We conclude with a discussion of future work in Section~\ref{sec:conclusion}.

\section{Examples} \label{sec:models}

Kinetic equations and elliptic equations with oscillating media are two examples of multiscale PDEs, for which computational schemes were developed separately.
The specific features of these problems were incorporated into the design of asymptotic preserving schemes and numerical homogenization methods, respectively.
We review these techniques and highlight the shared low-rank property of these two problems.

\subsection{Kinetic equations and asymptotic preserving methods}
Kinetic equations, which originate from statistical mechanics,  describe the evolution of probability density for identical particles in phase space.
A model equation, the radiative transfer equation (RTE), characterizes the evolution of photon density.
In the steady state, this equation is
\begin{equation}\label{eqn:rte}
-v\cdot\nabla_x u^\eps +\Ss^\epsilon[u^\eps] = f(x,v)\,, \quad (x,v)\in\Kc\times\Vb \,,
\end{equation}
where $f(x,v)$ is the light source, and the linear collision operator $\Ss^\eps$ describes the interaction of photons with the optical media.
The small parameter $\eps$ is encoded in this operator.

The operator  $\Ss^\eps$  defines several distinct regimes.
In the optically thick regime, it is defined by
\begin{equation}\label{eqn:collision}
\begin{aligned}
\Ss^\eps u(x,v) = \frac{1}{\eps}\int_\Vb & k(x,v,v')u(x,v')\rmd v' \\
                                        & - \frac{1}{\eps}\int_\Vb k(x,v',v) u(x,v) \rmd v' \,.
\end{aligned}
\end{equation}
In this case, $k(x,v,v')$ is the scattering coefficients that describes the possibility of a photon located at $x$ changing its velocity from $v$ to $v'$, and the parameter $\eps$ is called the {\em Knudsen number}, standing for the ratio of the mean free path to the typical domain length.
When the medium is optically thick, the mean free path is small, with $\eps \ll 1$.
This means the photon particles are scattered fairly often, and the system statistically achieves the equilibrium state, which can itself be characterized mathematically.
One example is to observe light in atmosphere, where the average mean free path is about $10$ m, and the observation is conducted at the scale of $10$ km, leading to $\eps\sim 10^{-3}$.
By performing asymptotic expansion in terms of $\eps$, the inhomogeneity in the velocity domain vanishes, and one can show that $u^\varepsilon(x,v)$ asymptotically approximates $u^\ast(x)$, a function without dependence on $v$ that solves a diffusion equation.
We have the following result from \cite{MR2839402}.
\begin{theorem}\label{thm:hom_transport}
Suppose that $u^\eps$ solves~\eqref{eqn:rte} with collision term $\Ss$ being isotropic, that is, $k(x,v,v')=\sigma(x)$ for some function $\sigma$.
Let $\mathcal{K}\subset \mathbb{R}^d$ be bounded with smooth boundary, and $\Vb = \Sb^{d-1}$.
Assume that the boundary condition is
\begin{equation}\label{eqn:transportBdy}
u^\eps(x,v) = \phi(x,v) \quad\text{on}\quad x\in\partial\mathcal{K}\,, \;\; v\cdot n_x<0\,.
\end{equation}
Then
\begin{equation}\label{eqn:conv_rte_diff}
\|u^\eps - u^\ast\|_{L_2(\rmd x \rmd v)} \to 0\,,
\end{equation}
where $u^\ast=u^\ast(x)$ solves
\begin{equation}\label{eqn:diff}
\nabla_x\cdot \left(\frac{1}{\sigma(x)}\nabla_x u^\ast(x) \right) =g(x)\,,\quad x\in\mathcal{K}\,,
\end{equation}
with the boundary condition
\begin{equation*}
u^\ast(x) = \xi_\phi(x)\,,\quad\text{on}\quad x\in\partial\mathcal{K}\,,
\end{equation*}
where $\xi_\phi(x)$ solves a proper  boundary layer equation and $g$ can be obtained from $\int f(x,v)\rmd{v}$.
\end{theorem}
This result indicates that the homogenized operator as $\eps \to 0$ is $\mathcal{L}^\ast \propto \nabla_x\cdot\left(({1}/{\sigma})\nabla_x\right)$. Similar results, when $k(x,v,v')$ fails to have the form of $\sigma(x)$ in the anisotropic optical media, are still available, but the explicit form of $\mathcal{L}^\ast$ is no longer available.

A second regime of interest for $\Ss^\eps$ is one in which the media is highly heterogeneous~\cite{MR1760042}:
\begin{equation}\label{eqn:collision2}
\begin{aligned}
\Ss^\eps u(x,v) = \int_\Vb & k\left(\frac{x}{\eps},v,v'\right)u(x,v')\rmd v' \\
                                        & - \int_\Vb k\left(\frac{x}{\eps},v',v\right) u(x,v) \rmd v' \,.
\end{aligned}
\end{equation}
In this case, the photons go through the media that oscillates at a small scale: For example, sunlight passing through heavy cloud with a large number of small droplets or laser beam passing through crystals. The amplitude of $k$ determines the photon scattering frequency. Since $k$ oscillates rapidly, photons also change rapidly between the high- and low-scattering  regimes. On a large scale, the photons can be viewed approximately as scattering with an averaged frequency. A mathematical result is as follows~\cite{MR1760042}.
\begin{theorem}\label{thm:hom_transport_2}
Let the conditions from Theorem~\ref{thm:hom_transport} hold, and suppose that the collision term $\Ss^\eps$ is defined in~\eqref{eqn:collision2}. Then
\begin{equation}\label{eqn:conv_rte_multi}
\|u^\eps - u^\ast\|_{L_2(\rmd x \rmd v)} \to 0\,,
\end{equation}
where $u^\ast(x,v)$ solves
\begin{equation}\label{eqn:hom_transport_2}
-v\cdot\nabla_x u^\ast +\Ss^\ast[u^\ast] = f(x,v)\,, \quad (x,v)\in\Kc\times\Vb \,,
\end{equation}
where $\Ss^\ast u(x,v) = \sigma^\ast(x)\int_\Vb u(x,v') - u(x,v) \rmd v'$ for some $\sigma^\ast(x)$. Furthermore, if $k(x,v,v')=\sigma(x)$ is periodic in $x$ with period $[0,1]^d$, then $\sigma^\ast=\int_{[0,1]^d} \sigma(x) \rmd x$.
\end{theorem}

In special cases, such as under periodic or random ergodic conditions, the function $\sigma^\ast$ can be computed explicitly.
(There are also works that investigate the asymptotic limit of the RTE when the system is both highly oscillatory and in diffusion regime; see \cite{MR1878799}.)

In both limiting regimes, the limiting equations~\eqref{eqn:diff} and~\eqref{eqn:hom_transport_2} can be solved much more efficiently than the original equation~\eqref{eqn:rte}.
The discretization of~\eqref{eqn:rte} is constrained strongly by $\eps$, due either to stability (as in~\eqref{eqn:collision}) or accuracy (as in~\eqref{eqn:collision2}).
By contrast, the solution $u^\ast$ varies smoothly, containing no $\eps$-scale effects, so can be obtained accurately by applying a discretization with mesh width $O(1)$ to the asymptotic limiting equation.
If the latter equation is available, computation of $u^\ast$ by this means is the recommended methodology.

Methods for kinetic equations are termed ``asymptotic preserving" (AP) if they can relax the  requirement that the discretization width $h$ satisfies $h =o(\eps)$ yet still capture the asymptotic limits.
Many different AP approaches have been proposed.
For linear equations, existing AP methods rely on even-odd or micro-macro decomposition.
For nonlinear equations, knowledge of the specific forms of the limits is usually required, and this knowledge is built into the solvers~\cite{MR3645390}.
As mentioned above, these specific forms are often not available, so many AP methods cannot be applied to a large set of multiscale kinetic equations.
This observation begs the question: Knowing the existence of the limit, but not its particular form, can we still devise efficient methods for solving kinetic equations?

\subsection{Elliptic equations and numerical homogenization}

Another class of multiscale equations that has been  investigated deeply is elliptic equations with highly oscillatory coefficient.
These problems have the form
\begin{equation}\label{eqn:elliptic}
-\nabla_x\cdot\left(a^\eps\left(x\right)\nabla_x u^\eps\right) = f \,,
\end{equation}
where $\eps \ll 1$ is the scale on which the media oscillates.
(The source term $f$ has no small-scale contribution.)

This equation is a model problem from petroleum engineering where it is crucial to precompute the underground flow before  expensive construction of infrastructure takes place~\cite{MR2801210}.
The problem is typically solved on  kilometer-scale domains, but the heterogeneities in the media can scale at centimeters.
Certain forms of this equation can be approximated effectively by an equation that can be solved efficiently.
Suppose the media coefficient $a^\eps(x)$ has the form $a(x,x/\eps)$, that is, it varies on two scales ($1$ and $\eps$), and moreover is periodic with respect to the fast variable (the second argument in $a(x,x/\eps)$).
Then in the limiting regime as $\eps\to0$, the solution $u^\eps$ converges to that of a homogenized equation, with the media ``smoothed-out,'' as described in the following result~\cite{MR1185639}.
\begin{theorem}\label{th:3}
Let $u^\eps$ solve~\eqref{eqn:elliptic} in the domain $x\in\Kc$ with zero boundary condition. Suppose $a(x,x/\eps)$ is periodic with respect to the second argument.
Then
\begin{equation}\label{eqn:conv_de}
\| u^\eps - u^\ast \|_{L^2} \lesssim \eps \left\|u^\ast\right\|_{H^2}\,,
\end{equation}
where $u^\ast$ solves the following effective equation with zero boundary condition:
\begin{equation}\label{eqn:elliptic_homo}
-\nabla_x\cdot(a^\ast(x)\nabla_x u^\ast) = f  \,, \quad  x\in\Kc \,,
\end{equation}
where $a^\ast$, the effective media, can be computed from a cell problem (See Definition 2.1 in \cite{MR1185639}).
\end{theorem}

As in the previous section, when a limiting equation can be derived explicitly, the best course for obtaining a useful solution it to solve this equation directly, as the mesh width in the  discretization scheme can be much larger than $\eps$. See \cite{MR2477579} for a discussion of a reduced number of basis functions and \cite{MR2830582} for computation of the effective media.

However, the validity and the specific form of the effective limit are known only in special cases like the one described in Theorem~\ref{th:3}.
In other cases, we seek a solver that relies on as little analytical knowledge as possible.
An approach known as {\em numerical homogenization} has been investigated extensively.
This approach is founded on two principles: a discretization scheme independent of $\eps$, and a numerical solution scheme that captures the true limiting behavior of the solution on the discrete level.
Variants of numerical homogenization include application of the $\Hc$-matrix, a purely algebraic technique \cite{MR3445676}; and a Bayesian approach that views the source $f$, and hence the solution $u^\eps$, as Gaussian fields~\cite{MR3369060}, which further translates to game theory \cite{MR3971243}. All these methods are successful, but they all implicitly rely on properties of the underlying elliptic equation. Can we devise an approach that applies to general problems with oscillatory media that exploits the low-rank property in the solution space, without using analytical structure explicitly?

\section{A unified framework for multiscale PDEs based on random sampling}\label{sec:framework}

We have given several examples of  multiscale models that arise in applications, and  mentioned several algorithmic approaches that make use of the limiting equations, when available.
We describe next the foundations of a unified scheme  that captures asymptotic limiting behavior automatically, even when the asymptotic limits are unavailable.
Our method exploits low-rank structure and uses random sampling to discover this structure.
We describe the low-rank property in Section~\ref{sec:rank} and the randomized SVD method for revealing this structure in Section~\ref{sec:rSVD}.

\subsection{Numerical rank}\label{sec:rank}

We consider a bounded linear operator $\mathcal{A}$, which maps $f\in\Xc$ to a space $\Yc$, that is
\[
\begin{aligned}
\Ac: & \quad \Xc &\rightarrow &\quad \Yc \\
& \quad f 				&\mapsto & \quad u.
\end{aligned}
\]
In the PDE setting, $\Ac$ is the solution operator that maps the boundary conditions and/or source term $f$ to the solution $u$.
The numerical rank of such an operator is defined  as follows.
\begin{definition}[Numerical rank]\label{def:num_rank}
	The numerical $\tau$-rank of $\Ac$ is the rank of the lowest-rank operator within the $\tau$-neighborhood of $\Ac$, that is,
	\begin{equation*}
	k_\tau(\Ac) := \min\{\dim \ran \tilde{\Ac}: \tilde{\Ac}\in\Lc(\Xc,\Yc), \|\tilde{\Ac}-\Ac\|\leq\tau  \}\,.
    \end{equation*}
	In other words, $k_\tau(\Ac)$ is this smallest dimension of the range among all the operators within distance $\tau$ of $\Ac$.
\end{definition}

When $\mathcal{A}$ is the PDE solution map, then $\tilde{\Ac}$ with low rank is also a linear map with a finite dimensionality.
It can be viewed as the discrete version (or a matrix of dimension $k_\tau(\Ac)$)
that approximates $\mathcal{A}$ within $\tau$ accuracy.
The definition suggests that if $\tilde{\Ac}$ can be found, it is optimal in the sense of numerical efficiency.
The concept is rather similar to the Kolmogorov $N$-width, defined as follows.
\begin{definition}[Kolmogorov $N$-width]\label{def:Kol_N}
	Given the linear operator $\Ac:\Xc\to\Yc$, the Kolmogorov $N$-width $d_N(\Ac)$ is the shortest distance from its range to all $N$-dimensional space, that is,
	\begin{equation}
    \begin{aligned}
	d_N(\Ac) :&= \min_{S:\dim S = N} d(\Ac,S)\\
    &=\min_{S:\dim S = N} \sup_f \min_{v \in S} \frac{\| \Ac f - v\|_\Yc }{\| f \|_{\Xc}} \,.
	\end{aligned}
    \end{equation}
\end{definition}
Indeed, the  Kolmogorov $N$-width  and numerical rank are related by the following result~\cite{MR4155236}.
\begin{proposition}
	For any linear operator $\Ac:\Xc\to\Yc$, we have the following.
	\begin{enumerate}[(a)]
		\item If the numerical $\tau$-rank is $N$, then $d_N(\Ac)\leq\tau$.
		\item If $d_{N}(\Ac)\leq \tau < d_{N-1}(\Ac)$, then the numerical $\tau$-rank is $N$.
	\end{enumerate}
\end{proposition}

For the three examples presented in Section~\ref{sec:models}, the numerical ranks can be calculated from their limiting equations.
For one-dimensional RTE in the diffusion regime, if we denote by $\Ac^\eps$ and $\Ac^\ast$ the solution operators of~\eqref{eqn:rte} and~\eqref{eqn:diff}, respectively, then noting that $\Ac^\ast$ can be approximated using $1/\sqrt{\tau}$ grid points to achieve $\tau$ accuracy, when $\eps<\tau$, the numerical rank is naturally $k_\tau(\Ac^\eps) \lesssim 1/\sqrt{\tau-\epsilon}$.
Without employing the knowledge of the existence of the limit, however, a brute-force discretization naturally requires $O(1/\eps\tau^{\alpha+1})$ degrees of freedom:  $O(1/\eps\tau)$ for the upwind discretization in $x$ and $O(1/\tau^\alpha)$ for the discretization in $v$, where $\alpha$ depends on the particular numerical integral accuracy.
Translating into Green's-matrix language, this observation means that $\mathsf{G}^\eps$ is represented by $O(1/\eps\tau^{\alpha+1})$ degrees of freedom but its range can be captured by a compressed Green's matrix $\Gs^\ast$ with just $O(1)$ column vectors.

The same argument applies to the elliptic equation on a two-dimension domain with high oscillations.
When second-order linear finite elements are used, with no knowledge of the limiting system, $O(1/\eps^2\tau)$ degrees of freedom are required, dropping to $O(1/\tau)$ when the the limiting system is known.
In other words, the full Green's matrix $\Gs^\eps$ requiring $O(1/\eps^2)$ degrees of freedom can be well-represented using just $O(1)$ column vectors.

In all these cases, the degrees of freedom for a given numerical method are substantially larger than the numerical rank of the problem.
Thus, much of the information in these full-blown representations is redundant and compressible.
A low rank representation exists and yields a much more economical representation.

\subsection{Random sampling in numerical linear algebra} \label{sec:rSVD}

Knowing the existence of the low rank structure and finding such a structure are very different goals.
The Kolmogorov $N$-width is a concept developed in but has made little impact in numerical PDEs for a  simple reason: Traditional PDE solvers require a predetermined set of basis functions, while the Kolmogorov $N$-width looks for ``optimal'' basis functions.
How can an optimal  basis be found without first forming the full basis?
Translated to linear algebra, this question is about finding the dominant singular  vectors in a matrix without forming the whole matrix.
Specifically, if $\As\in\Rb^{m\times n}$ is known to be approximately low rank, meaning that there exists $\Us_r$, a $m \times r$ matrix with orthonormal columns with $r \ll \min(m,n)$ and
\begin{equation*}
\|\As-\As_r\| = \|\As - \Us_r\Us_r^\top\As\| \ll \|\As\|\,,
\end{equation*}
can we find $\Us_r$ without forming the full matrix $\As$?

In linear algebra, it is well-known that $\Us_r$ is simply the collection of the first $r$ singular vectors of $\As$. Writing
\begin{equation}
\label{eqn:svd}
\As = \Us\Sigma \Vs^\top = \sum_{i=1}^n \sigma_i u_i v^\top_i\,,
\end{equation}
where $\Us=\left[u_1\,,u_2\,,\dots,u_n\right]\in\Rb^{m\times n}$ and $\Vs=\left[v_1\,,v_2\,,\dots,v_n \right]\in \Rb^{n\times n}$ contain the left/right singular vectors and $\Sigma = \diag(\sigma_1,\sigma_2,\dotsc,\sigma_n)$ contains the singular values, then $\Us_r$ is the first $r$ columns in $\Us$.

The standard method for computing the SVD requires $\As$ to be stored and computed with. But the celebrated randomized SVD (rSVD) method~\cite{MR2806637} captures the range of a given matrix by means of random sampling of its column space, which requires only computation of matrix-vector products involving $\As$ and random vectors --- operations that can be performed without full storage or knowledge of $\As$.
Implementation of rSVD is easy and its performance is robust.

The idea behind the algorithm is simple. If matrix $\As\in\Rb^{m\times n}$ has approximate low rank $r\ll \min\{m,n\}$, the matrix maps an $n$-dimensional sphere to an $m$-dimensional ellipsoid that is ``thin:'' $r$ of its axes are significantly larger than the rest.
With high probability, vectors that are randomly sampled on the $n$-dimensional sphere are mapped to vectors that lie mostly in a $r$-dimensional subspace of $\Rb^m$ --- the range of $\As$.
An approximation to $\As_r$ can be obtained by projecting onto this subspace.

A precise statement of the performance of randomized SVD is as follows~\cite{MR2806637}.
\begin{theorem}
	\label{thm:average_spectral}
	Let $\As$ be an $m\times n$ matrix. Define
	\begin{equation}\label{eqn:define_Y}
	\Ys = \As\Omega\,,
	\end{equation}
	where $\Omega=\left[\omega_1\,,\dotsc,\omega_{r+p}\right]$ is a matrix of size $n\times (r+p)$ with its entries randomly drawn from an i.i.d.~normal distribution, where $p$ is an oversampling parameter. If $\sigma_{r+1}\ll\sigma_1=O(1)$, then the projection of $\As$ onto the space spanned by $\Ys$, defined by
\begin{equation*}
\Ps_\Ys(\As) = \Ys(\Ys^\top\Ys)^{-1}\Ys^\top\As\,,
\end{equation*}
yields that $\|\As - \Ps_\Ys(\As)\| \ll \sigma_1$ with high probability, and
\begin{equation*}
\mathbb{E}\,\|\As - \Ps_\Ys(\As)\| \lesssim \frac{r}{p-1}\sigma_{r+1} \ll\sigma_1 \,.
\end{equation*}
\end{theorem}

The result reconstructs the range of $\As$ in a nearly optimal way. It is optimal in efficiency because to capture a rank-$r$ matrix, only $r+p$ matrix-vector products involving $\As$ are required for the calculation of $\Ys$, and the oversampling parameter $p$ is typically quite modest.
($p=5$ is a typical value.)
The result is nearly optimal in accuracy as well.
The error bound relies only on $\sigma_{r+1}$, which is expected to be smaller than $\sigma_1$.
The decay profile of singular values do not affect the approximation accuracy.

If a low-rank approximation to  the matrix $\As$ is required, and not just an approximation of its range, another step involving multiplications with its transpose is needed.
The full method is shown in Algorithm~\ref{alg:rsvdMatrix}.

\begin{algorithm}[ht]
	\caption{Randomized SVD}\label{alg:rsvdMatrix}
	\begin{algorithmic}[1]
		\State Given an $m\times n$ matrix $\As$, target rank $r$ and oversampling parameter $p$;
		\State Set $k = r+p$;
		\State \textbf{Stage A:}
        \Indent
		\State Generate an $n\times k$ Gaussian test matrix $\Omega$;
		\State Form $\Ys =  \As \Omega$;
		\State Perform the QR-decomposition of $\Ys$: $\Ys=\Qs\Rs$
        \EndIndent
		\State \textbf{Stage B:}
        \Indent
		\State Form $\Bs =  \As^\top \Qs$;
		\State Compute the SVD of the $k \times n$ matrix $\Bs^\top = \widetilde{\Us}\Sigma \Vs^\top$;
		\State Set $\Us = \Qs \widetilde{\Us}$;
        \EndIndent
		\State \textbf{Return:} $\Us, \Sigma, \Vs$.
	\end{algorithmic}
\end{algorithm}

\section{Random sampling for multiscale computation}

Here we describe how rSVD can be incorporated into multiscale PDE solvers to exploit the low-rank structure of these equations.
Our procedure is composed of both offline and online stages.
Low rank structure is learned in the offline stage, while in the online stage, the solution for the given source / boundary term $f$ in~\eqref{eqn:general} is extracted.

We consider in particular the following boundary value problem:
\begin{equation} \label{eq:bvp}
\begin{cases}
(\Lc^\eps u^\eps)(x) = 0\,, \quad x\in \Kc \,, \\
\Bc u(x) = \phi(x) \,, \quad x\in \partial\Kc \,,
\end{cases}
\end{equation}
where $\Bc$ is the boundary condition operator, $\partial\Kc$ the boundary associated with domain $\Kc$, and we now denote the source term (the boundary data) by $\phi$. Our fundamental  goal is to construct the low-rank approximation to the Green's operator $\Gs^\eps$ for \eqref{eq:bvp}. With this operator in hand, the solution can be computed for any value of the boundary conditions $\phi$ at a relatively small incremental cost.

If we apply rSVD to approximate $\Gs^\eps$ directly, we need to compute products of this operator with random vectors.
For the problem \eqref{eq:bvp}, this operation corresponds to solving the problem with $\phi$ replaced by random boundary conditions.
Even to solve one such problem efficiently is a computationally challenging task.
We use the domain decomposition framework.

\begin{figure}
	\centering
	\includegraphics[width=0.4\textwidth]{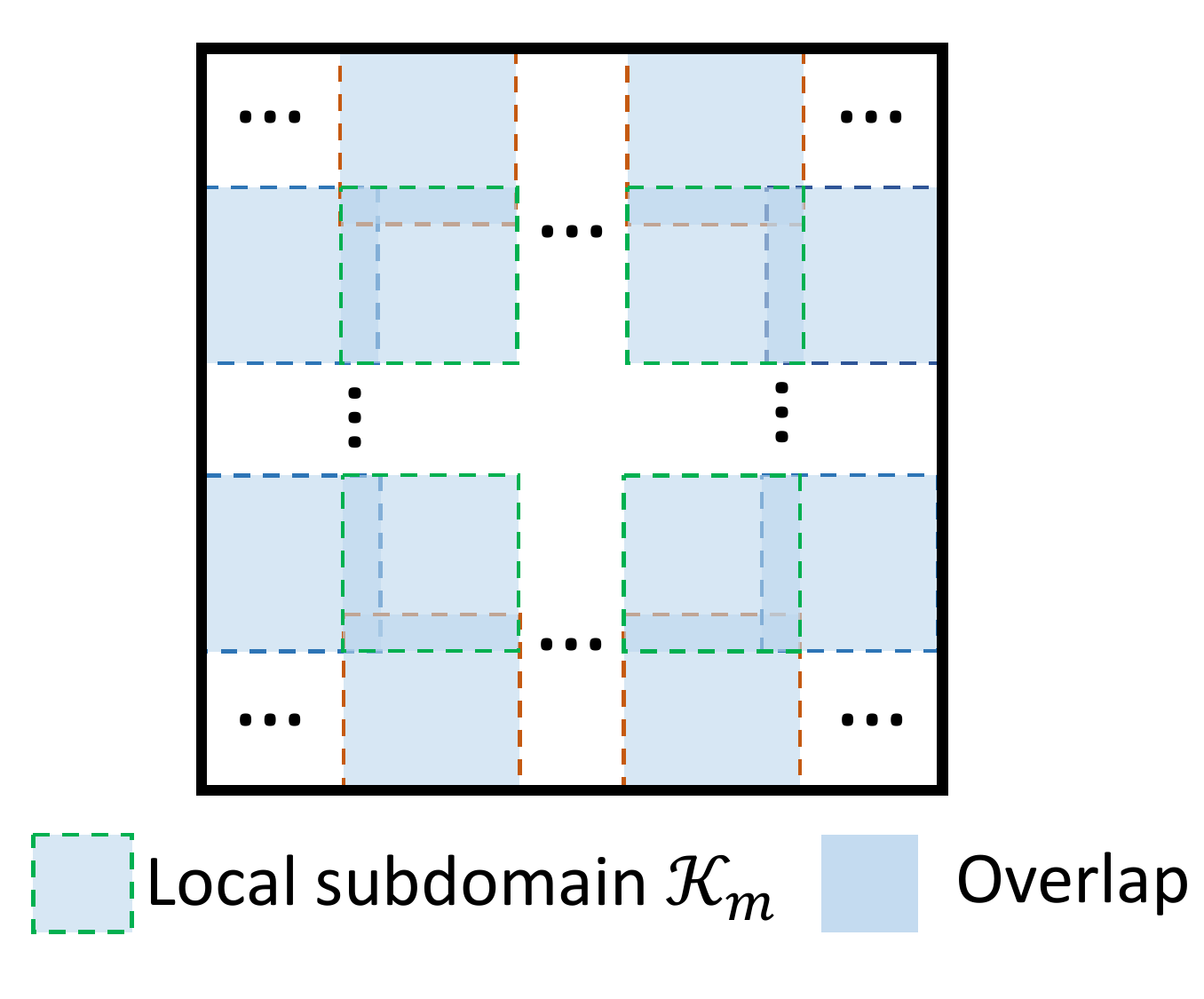}
	\caption{Illustration of domain decomposition of a rectangular domain $\Kc$ in 2D with overlap.}
	\label{fig:dd2}
\end{figure}

We start by partitioning the domain $\Kc$ into subdomains as follows:
\begin{equation}\label{eqn:partition}
\Kc = \bigcup_{m=1}^M \Kc_m\,,
\end{equation}
where the patches $\Kc_m$ overlap, in general.
We denote by $\partial \Kc_m$ the boundary associated with $\Kc_{m}$.
Furthermore, we identify the subregions that intersect with $\Kc_m$ as follows:
\[
\Ic_m = \{n\in\Nb: 1\leq n \leq M\,, \Kc_m\cap\Kc_n \neq \varnothing\}\,,
\]
and define the interior of the patch to be
\[
\widetilde{\Kc}_m = \Kc_m \backslash \left(\bigcup_{n\in\Ic_m}\Kc_n\right)\,.
\]
For this particular partition of the domain, we define the partition-of-unity functions $\chi_m$, $m=1,2,\dotsc,M$ to have the following properties:
\begin{equation}\label{eqn:POU}
\begin{aligned}
&\sum_{m = 1}^M \chi_m(x)= 1\,, \quad \forall x\in\Kc\,,\\
&\text{with}\quad
\begin{cases}
0 \leq \chi_m(x) \leq 1\,,  & x\in\Kc_m\\
\chi_m(x) = 0\,,  & x\in\Kc\backslash\Kc_m.
\end{cases}
\end{aligned}
\end{equation}
We choose a discretization that resolves the small scales in the solution, defining a mesh width $h\ll \eps$.
(The number of subdomains $M$ is independent of $\eps$.)
A typical decomposition is illustrated in Figure~\ref{fig:dd2}.

How do we design the offline stage to ``learn'' the low-rank approximation?
We propose two approaches that lead to two different kinds of algorithms.
In the first approach we learn the optimal basis functions within each subdomain, while the second algorithm employs Schwarz iteration, preparing the boundary-to-boundary map in the offline stage.
Other PDE solvers that utilize randomness can also be found in~\cites{MR3477310,MR3824169,MR4050504}. In particular in~\cite{MR3477310} the authors studied, specifically for elliptic type equations, the generalized eigenvalue problem of the stiffness and mass matrices, and give an error bound using the largest eigenvalue obtained offline.

\subsection{Learning basis functions}\label{sec:learning_basis}

In standard domain decomposition, the local discretized Green's matrix $\Gs_m$ is assembled from a ``full"  collection of basis functions in the patch $\Kc_m$.
The global solution to \eqref{eq:bvp}, confined to each $\Kc_m$, is a linear combination of the columns $\Gs_m$.
The coefficients of these combinations are chosen so that that the continuity conditions across patches and the exterior boundary condition are all satisfied.
The complete process can be outlined as follows.
\begin{enumerate}[(1)]
	\item \textbf{Offline stage}: For $m=1,2,\dots,M$, find
	\begin{equation*}
	\Gs_m = \left[b_{m,1}\,,b_{m,2}\dots \right]\,,
	\end{equation*}
	where each local function $b_{m,n}$ is a solution to \eqref{eq:bvp} restricted to the subdomain $\Kc_m$, with fine grid $h \ll \eps$ and delta-function boundary conditions. That is,
	\begin{equation}\label{eqn:full_basis}
    \begin{cases}
    \Lc^\eps b_{m,n} = 0\,, \quad & x\in \Kc_m \\
    b_{m,n} = \delta_{m,n}\,,\quad  & x\in \partial\Kc_m\,,
    \end{cases}
    \end{equation}
    where $\delta_{m,n}$ is the Kronecker delta that singles out the $n$-th grid point on the boundary $\partial \Kc_m$.
\item \textbf{Online stage}: The global solution is
	\[
	u = \sum_{m=1}^M u_m\chi_m = \sum_{m=1}^M \chi_m\Gs_{m}{c}_m,
	\]
	with the support of each $u_m=\Gs_{m}{c}_m$ confined to $\Kc_m$, where $c_m$ is a vector of coefficients determined by the boundary conditions $\phi$ and continuity conditions across the patches.
\end{enumerate}

The complete basis represented by $\Gs_m$ has a low-rank structure that can be revealed using randomized SVD.
Instead of using delta functions as the boundary conditions, we propose to obtain basis functions by setting random values on $\partial \Kc_m$, as follows:
\begin{equation}\label{eqn:rand_basis}
\begin{cases}
\Lc^\eps r_{m,n} = 0\,, \quad & x\in \Kc_m, \\
r_{m,n} = \omega_{m,n}\,,\quad & x\in \partial\Kc_m,
\end{cases}
\end{equation}
where $\omega_{m,n}$ is defined to have a random value drawn i.i.d. from a normal distribution at each grid point in $\partial\Kc_m$.
Denoting $\Gs_m^\rmr=\{r_{m,1}\,,r_{m,2}\,,\cdots\}$, we have from linearity of the equation that
\begin{equation*}
\Gs_m^\rmr = \Gs_m\Omega\,,
\end{equation*}
where $\Omega$ is a random i.i.d. matrix with entries $\omega_{m,n}$.
This $\Gs_m^\rmr$ is used in the online stage, as an accurate surrogate of $\Gs_m$, see Algorithm~\ref{alg:summary}.

\begin{algorithm}
	\caption{A general framework for multiscale PDE $\Lc^\eps
		u^\eps = 0$ over $\Kc$ with $\Bc u^\eps = f$ on $\partial\Kc$}\label{alg:summary}
	\begin{algorithmic}[1]
		\State \textbf{Domain Decomposition}
		\Indent
		\State Partition domain according to~\eqref{eqn:partition}.
		\EndIndent
		\State \textbf{Offline Stage:}
		\Indent
		\State Prepare i.i.d. Gaussian vectors $\omega_{m,i},i=1,\ldots,k_m$ on each $\partial\Kc_m$.
        \State Solve the basis function $r_{m,i}$ in~\eqref{eqn:rand_basis} on each $\Omega_m$, and collect the local basis in $\Gs_m$.
		\EndIndent
		\State \textbf{Online Stage:}
		\Indent
		\State Use continuity condition and global boundary data $\phi$ to determine coefficient vectors $c_1,c_2,\dotsc,c_M$, and set
        \begin{equation}\label{eqn:online}
        u = \sum_{m=1}^{M}\chi_m\Gs_{m} c_m
        \end{equation}
		\EndIndent
		\State \textbf{Return:} approximate global solution $u$.
	\end{algorithmic}
\end{algorithm}

\begin{figure}[tb]
	\centering
	\includegraphics[width = 0.23\textwidth]{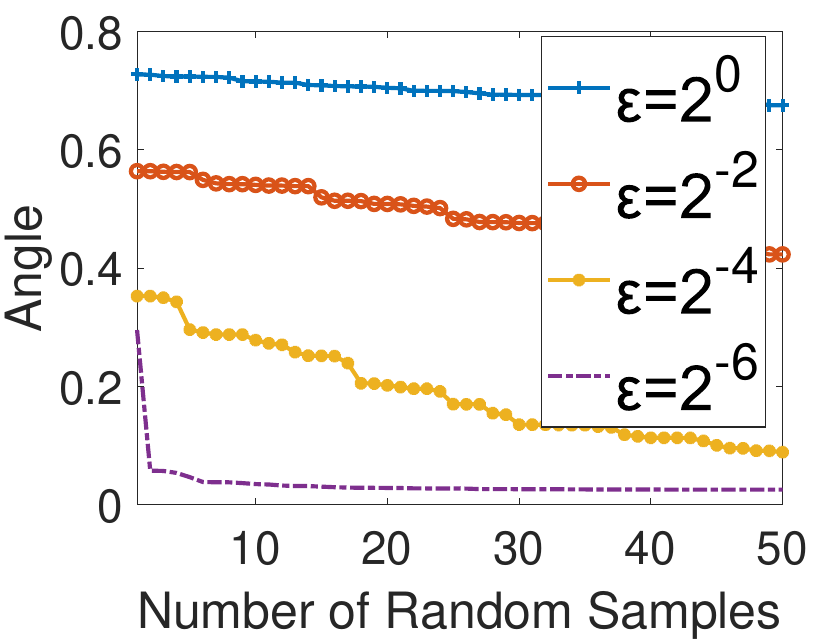}
    \includegraphics[width = 0.23\textwidth]{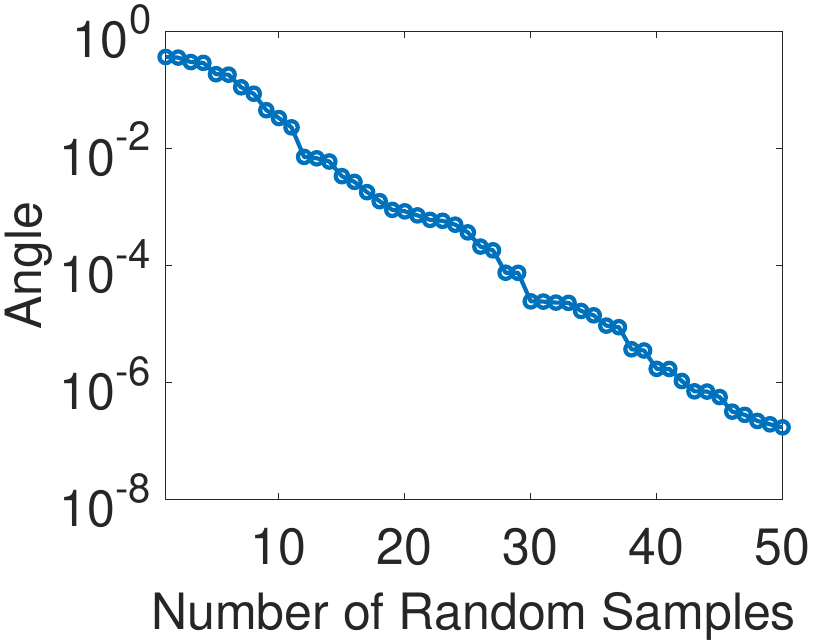}
	\caption{Angle between the true Green's matrix $\Gs_m$ and the approximate version $\Gs_m^\rmr$, confined on the interior of a patch $\widetilde{\Kc}_m$ for some $m$, as the number of random samples increases. Left plot: Angle for 1D RTE~\eqref{eqn:rte} with diffusive kernel~\eqref{eqn:collision} and various values of $\eps$. Right plot: Angle for elliptic equation~\eqref{eqn:elliptic} on a rectangular local patch with $\eps = 2^{-4}$.}
	\label{fig:range_local}
\end{figure}

Although we do not apply full-blown rSVD here, the homogenizable and low-rank property of the local solution space implies that  $\Gs_m^\rmr$ and $\Gs_m$ share similar range with the number of basis functions $k_m$ in $\Gs_m^\rmr$ being much smaller than $n_m$, the number of grid points on $\partial \Kc_m$, and independent of $\eps$.
In Figure~\ref{fig:range_local} we plot the angles between $\Gs_m^\rmr$ and $\Gs_m$ for two of the equations discussed in Section~\ref{sec:models}. In both cases, and for small $\eps$, the approximated Green's matrix quickly recovers the true Green's matrix as the number of samples $k_m$ increases, and thus captures the local solution space. We should note that if $\Gs_m$ does not have low-rank structure, in the sense that $k_m\sim n_m$, then solving~\eqref{eqn:rand_basis} would be equally expensive as solving~\eqref{eqn:full_basis}, hence the random sampling technique does not gain any computational efficiency when the system is not homogenizable.

\begin{figure}[tb]
	\centering
	\includegraphics[width = 0.23\textwidth]{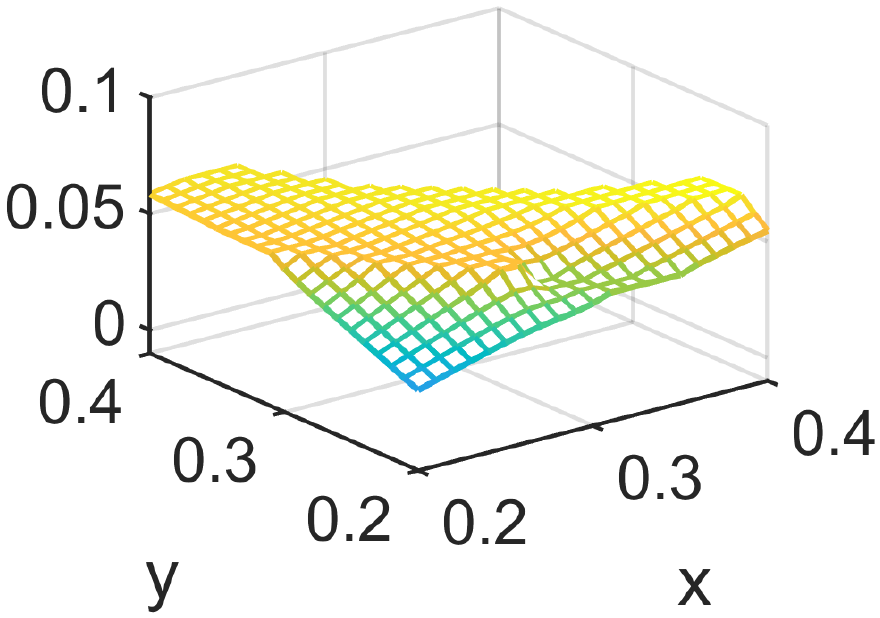}
    \includegraphics[width = 0.23\textwidth]{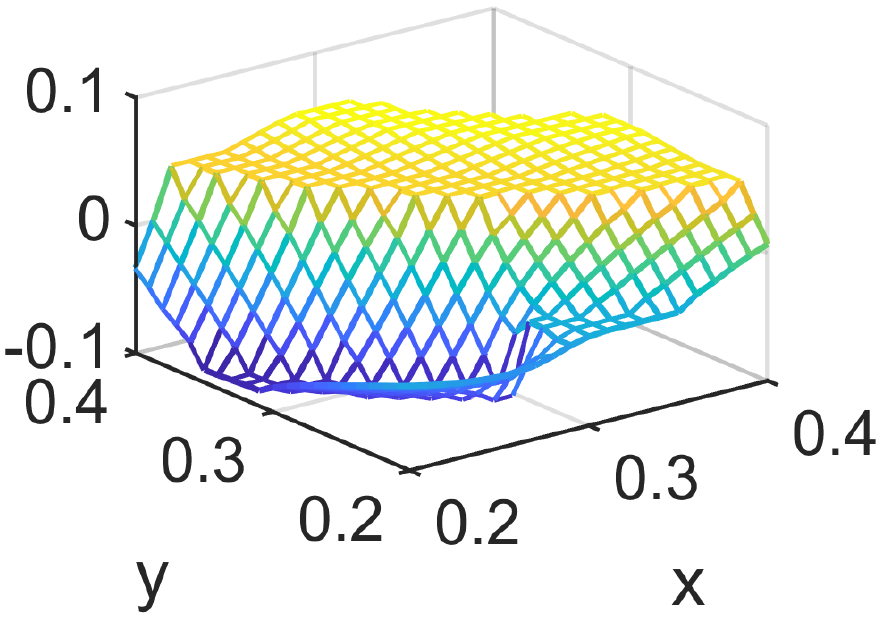}
    \includegraphics[width = 0.23\textwidth]{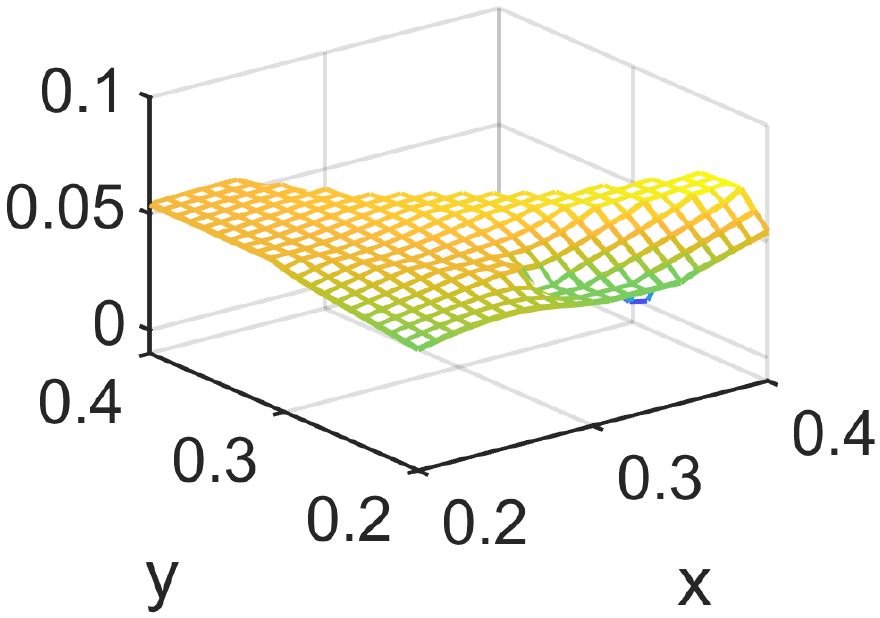}
    \includegraphics[width = 0.23\textwidth]{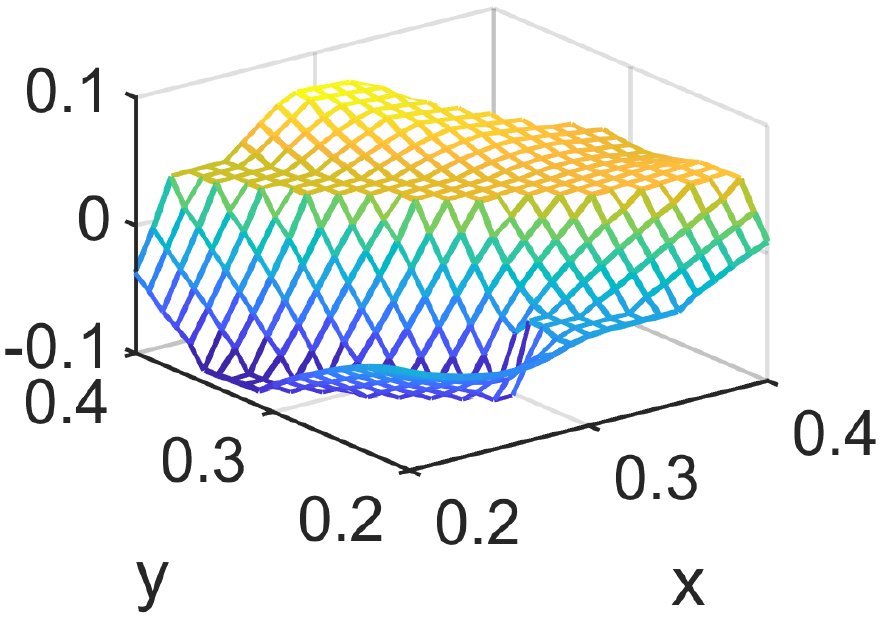}
	\caption{Optimal basis functions and their projections onto the approximate spaces. First row: First two singular vectors of the Green's matrix on a local patch. Second row: Projections onto the space spanned by 6 random sampled basis functions. Note that the small random sample captures well the leading eigenvectors of the true Green's operator.}
	\label{fig:local_basis}
\end{figure}

In Figure~\ref{fig:local_basis} we showcase the basis functions on a patch for elliptic equation with media coefficient $a(x_1,x_2) =1 + 1000 \, \mathbf{1}_S(x_1,x_2)$, with $S = \{(x_1,x_2)\in [0,1]^2:(x_1\cos(100\sqrt{(x_1-0.5)^2+(x_2-0.5)^2}))\leq x_2-0.5 \}$. For this non-conventional media without any periodic structure, the traditional multiscale methods are no longer valid, but our method still quickly captures the optimal basis.

\begin{figure}[tb]
	\centering
	\includegraphics[width = 0.23\textwidth]{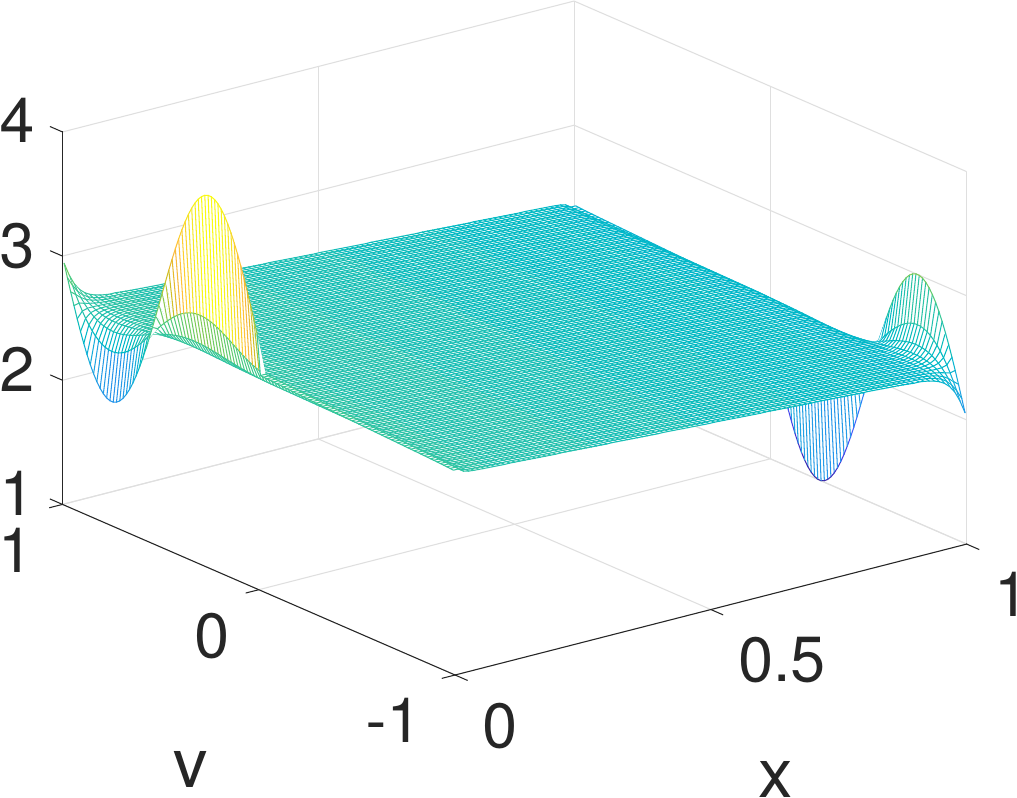}
    \includegraphics[width = 0.23\textwidth]{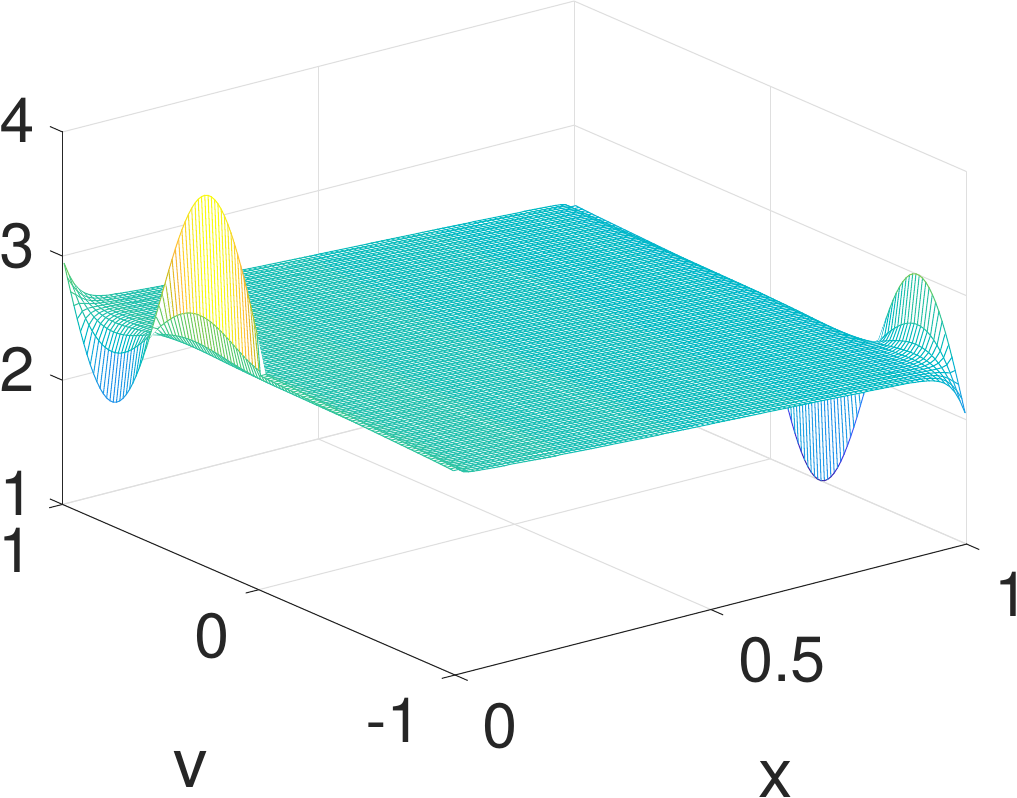}
    \includegraphics[width = 0.23\textwidth]{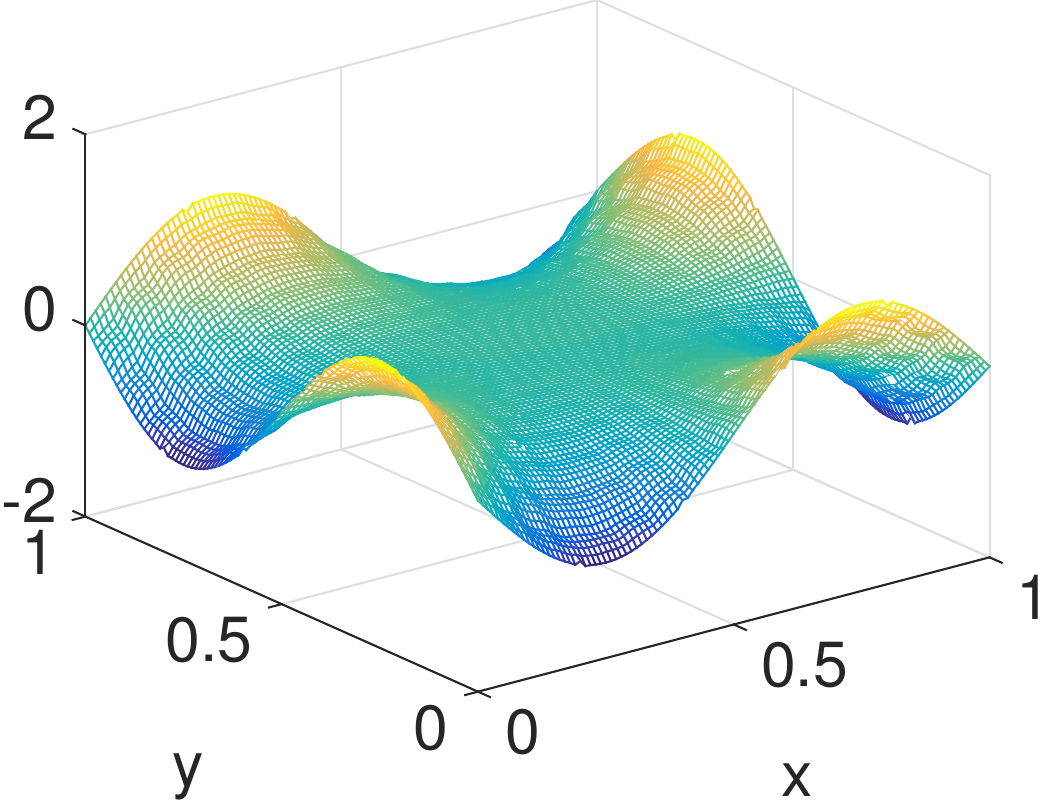}
    \includegraphics[width = 0.23\textwidth]{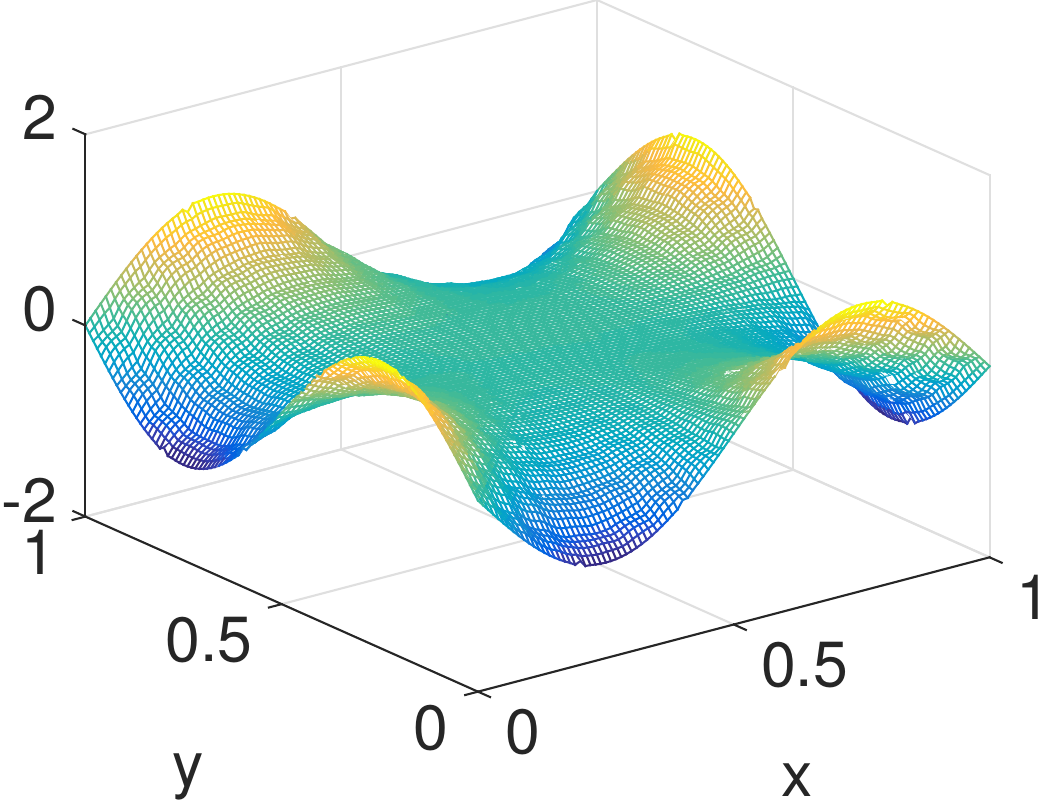}
	\caption{First row: Solutions for 1D RTE~\eqref{eqn:rte} with Henyey-Greenstein kernel in~\eqref{eqn:collision} with $\eps=2^{-6}$. Second row: Solutions for the elliptic equation~\eqref{eqn:elliptic} with Dirichlet boundary condition and highly oscillatory medium
$a(x,x/\eps)= 2 + \sin(2\pi x_1) \cos(2\pi x_2) +
\frac{2 + 1.8\sin(2\pi x_1/\eps)}{2 + 1.8\cos(2\pi x_2 / \eps)} +
\frac{2 + \sin(2\pi x_2/\eps)}{2 + 1.8\cos(2\pi x_1/\eps)}$ with $\eps = 2^{-4}$. The left column shows reference solutions while the right column is obtained from randomized reduced bases.}
	\label{fig:global}
\end{figure}

For particular boundary conditions $\phi$, the global solution is assembled from the local basis functions in the online stage.
Two numerical examples are shown in Figure~\ref{fig:global}.
In both examples, there is little visible difference between the reference solution and the approximated one computed from the reduced random basis.
Only 8.3\% and 62.5\%, respectively, of the degrees of freedom required by the full basis are needed to represent these solutions using a random basis; see details in~\cite{MR4155236}.

Since we do not have access to the full set of basis functions, the condition that $u$ defined in~\eqref{eqn:online} is continuous across subdomains can be satisfied only in a least-squares sense; see \cite{MR4155236} for details.

\subsection{A low-rank Schwarz method}\label{sec:Schwarz}

Our second approach for exploiting the low-rank property in multiscale computations is based on Schwarz iteration. The Schwarz method is a standard iteration algorithm within the domain decomposition framework, in which boundary-value problems are solved on the patches, with neighboring patches subsequently exchanging information and re-solving until consistency is attained. The exchange of boundary information between neighboring patches is known as the boundary-to-boundary (BtB) map. The map has an exploitable low-rank property.

To develop the approach, we write the solution of \eqref{eq:bvp} as
\begin{equation}\label{eqn:global}
u^\eps(x) = \sum_{m=1}^N \chi_m(x) u^\eps_m(x)\,,
\end{equation}
where the partition-of-unity functions $\chi_m$ are defined in~\eqref{eqn:POU}. The solution $u^\eps_m$ on patch $m$ is uniquely determined by $f_m$, its local boundary condition, according to the equation
\begin{equation}\label{eqn:elliptic_local}
\begin{cases}
\mathcal{L}^\eps u^\eps_m= 0\,, & \quad x\in \Kc_m \\
u_m^\eps(x) = f_m(x) \,, & \quad x \in \partial\Kc_m\,.
\end{cases}
\end{equation}
The Schwarz method starts by initial guesses to the local boundary conditions $f_m=f_m^0$, then on iteration $t$, it solves the subproblems~\eqref{eqn:elliptic_local} with $f_m=f_m^t$ to obtain all local solutions $u_m^\eps$. By confining $u_m^\eps$ on the  boundaries of adjacent patches, one updates the boundary conditions for surrounding patches:
\begin{equation}\label{eqn:update2}
f^t_m(x) \xrightarrow{\Sc_m} u_m^\eps(x)|_{\widetilde{\Kc}_m} \xrightarrow{\Pc_m} f^{t+1}_n(x)\,, \forall n\in\Ic_m\,.
\end{equation}
Here $\Sc_m$ denote the solution to~\eqref{eqn:elliptic_local} confined in the interior of $\Kc_m$, and $\Pc_m$ takes the trace of the solution on the neighboring boundaries $\Kc_m\cap\partial\Kc_n$ for $n\in\Ic_m$, for the updated boundary condition.

Define the BtB map by $\Ac_m := \Pc_m \circ \Sc_m$, and define $\Ac$ and $f^t$ to be the aggregation of $\Ac_m$ and $f_m^t$, respectively, over $i=1,2,\dotsc,M$. We can then write the updating procedure as
\[
f^{t+1} = \Ac f^t\,.
\]
The overall method is summarized in Algorithm~\ref{alg:schwarz}.

\begin{algorithm}[ht]
	\caption{Schwarz method for multiscale PDE $\Lc^\eps u^\eps = 0$ over $\Kc$ with $\Bc u^\eps = f$ on $\partial\Kc$}\label{alg:schwarz}
	\begin{algorithmic}[1]
		\State Given total iterations $T$;
        \State \textbf{Domain Decomposition}
        \Indent
            \State Partition domain according to~\eqref{eqn:partition}.
        \EndIndent
        \State \textbf{Schwarz Iteration:}
        \Indent
		  \State Initialize $f_m^0$ for each $\partial\Kc_m$ and set $t = 0$\,.
		  \State \textbf{While} $|f_n^{t}-f_n^{t-1}|>\text{TOL}$
          \Indent
		      \State Solve~\eqref{eqn:elliptic_local} for $u_m^t$ using $f_m^t$ for each $m$\,.
		      \State Update $f_{n}^{t+1} = u_m^t$ on $\Kc_m \cap \partial\Kc_n$, $n\in\Ic_m$\,.
		      \State $t\to t+1$.
		  \EndIndent
		  \State \textbf{End}
		  \State Solve~\eqref{eqn:elliptic_local} for $u_m^t$ using $f_m^t$ for each $m$\,.
		  \State Assemble global solution $u=\sum_{m=1}^N\chi_m u_m^t$\,.
        \EndIndent
		\State \textbf{Return:} approximated global solution $u^T$.
	\end{algorithmic}
\end{algorithm}

Most of the computation in the Schwarz method during the iteration comes from solving the boundary-value PDEs on the patches, to implement the map $\Ac$.
Since the PDE is homogenizable, the solution space on each patch is approximately low rank, and the map $\Ac$ can be expected to inherit this property. If we can ``learn" this operator in an offline stage, and simply apply a low-rank approximation repeatedly in the online stage, the online part of Algorithm~\ref{alg:schwarz} can be made much more efficient.
In our approach, Algorithm~\ref{alg:rsvdMatrix} is used to compress the map $\Ac$.

This approach is quite different from the one described in Section~\ref{sec:learning_basis}, in the sense that it is not only the range of the solution space, but the whole operator that is being approximated.
To apply Algorithm~\ref{alg:rsvdMatrix}, we need to define the ``adjoint operator" for $\Ac$ on the PDE level. This operator is composed of the adjoints $\Sc_m^\ast$ for  the  local solution operators $\Sc_m$ of \eqref{eqn:elliptic_local} on each domain $\Kc_m$. The form of $\Sc_m^\ast$ is specific to the PDE; we use the elliptic equation as an example.
Defining $\mathcal{L}^\eps = \nabla\cdot\left(a(x,x/\eps)\nabla\right)$, $\Sc_m^\ast$ is defined in the following result.
\begin{theorem}\label{thm:ad}
	Let $\Sc_m$ be the confined solution operator for the elliptic equation with Dirichlet boundary condition on patch $\Kc_m$. Given any function $g$ supported on $\widetilde{\Kc}_m$, the adjoint operator $\Sc_m^\ast$ acting on $g$ is given by:
	\begin{equation}\label{eqn:thm3}
	{\Sc}_m^\ast g = a \frac{\partial h}{\partial n}\Big|_{\partial\Kc_m}\,,
	\end{equation}
	where $\frac{\partial h}{\partial n}$ is the outer normal derivative on $\partial \Kc_m$ and $h$ solves the following sourced elliptic equation:
	\begin{equation}\label{eqn:thm4}
	\begin{cases}
	\nabla\cdot \left(a\left(x,\frac{x}{\eps}\right)\nabla h(x)\right) = {g}\,, & x\in\Kc_m \\
	h(x) = 0\,, & x\in \partial\Kc_m\,.
	\end{cases}
	\end{equation}
\end{theorem}

We also describe calculation of the adjoint operator $\Sc_m^\ast$ for the RTE~\eqref{eqn:rte}.
\begin{theorem}\label{thm:adjoint_S}
Let $\Sc_m$ be the confined solution operator for RTE~\eqref{eqn:rte} and the conditions in Theorem~\ref{thm:hom_transport} hold. Given any function $g$ supported on $\widetilde{\Kc}_m\times\Vb$, the adjoint operator $\Sc_m^\ast$ is defined as follows
\begin{equation}\label{eqn:map_ad}
\Sc_m^\ast g(x,v) = h(x,v) \,, \quad x\in\partial\Kc_m\,, \quad v\cdot n_x<0 \,,
\end{equation}
where $h$ solves the adjoint RTE over $\Kc_m$, which is
\begin{equation}\label{eqn:RTE_ad}
-v\cdot \nabla_x h - \Sc^\eps[h] = g(x,v) \,, \quad (x,v)\in\Kc_m\times\Vb \,,
\end{equation}
with outgoing boundary condition $h(x,v) = 0$ on $x\in\partial\Kc_m$ and $v\cdot n_x>0$.
\end{theorem}

The specific form of the adjoint operator $\Sc_m^\ast$ allows us to adapt the randomized SVD algorithm to compress the confined solution map $\Sc_m$; see Algorithm~\ref{alg:rsvd}. This method requires only $k$ solves of local PDE~\eqref{eqn:elliptic_local} and sourced adjoint PDE~\eqref{eqn:thm4} (or~\eqref{eqn:RTE_ad}), together with a QR factorization and SVD of relatively small matrices. The overall low-rank Schwarz iteration is then summarized in Algorithm~\ref{alg:schwarz_red}.

\begin{algorithm}[ht]
	\caption{Randomized SVD for $\Sc_m$}\label{alg:rsvd}
	\begin{algorithmic}[1]
		\State Given target rank $r$ and oversampling parameter $p$;
		\State Set $k = r+p$;
        \State \textbf{Stage A:}
        \Indent
		  \State Generate $k$ random boundary conditions $\xi_j$ on $\partial\Kc_m$.
		  \State Solve~\eqref{eqn:elliptic_local} using $\xi_j$ as boundary conditions and restrict the solution over $\widetilde{\Kc}_m$ to obtain $u_j^\eps$.
		  \State Find orthonormal basis $\Qs = [q_1,\dots,q_k]$ of $\widetilde{U}=[u_1,\dots,u_k]$.
        \EndIndent
        \State \textbf{Stage B:}
		\Indent
		  \State Construct zero extension of $q_k$ over $\Kc_m$, denoted by $\widetilde{q}_k$.
		  \State Solve~\eqref{eqn:thm4} or~\eqref{eqn:RTE_ad} for $h_j$ using $\widetilde{q}_k$ as source.
		  \State Compute $b_j$ using $h_j$ by flux~\eqref{eqn:thm3} or restriction on incoming boundary~\eqref{eqn:map_ad}.
		  \State Assemble all fluxes $B = [b_1,\dots,b_k]$.
		  \State Compute SVD of $B^\ast = \widetilde{\Us}_k\Sigma_k \Vs_k^\ast$.
		  \State Compute $\Us_k = \Qs \widetilde{\Us}_k$.
        \EndIndent
		\State \textbf{Return} $\Us_k,\Sigma_k,\Vs_k$.
	\end{algorithmic}
\end{algorithm}

\begin{figure}[htb]
	\centering
	\includegraphics[width = 0.23\textwidth]{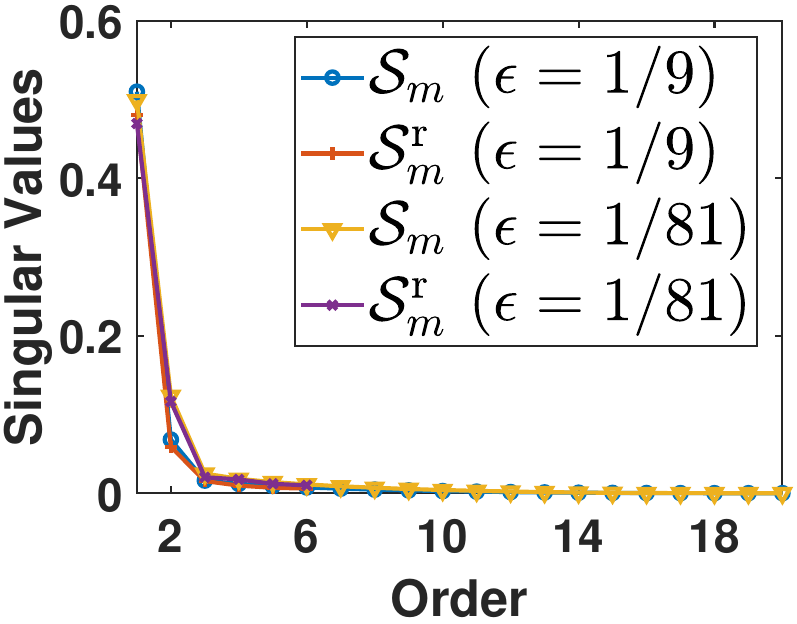}
	\includegraphics[width = 0.23\textwidth]{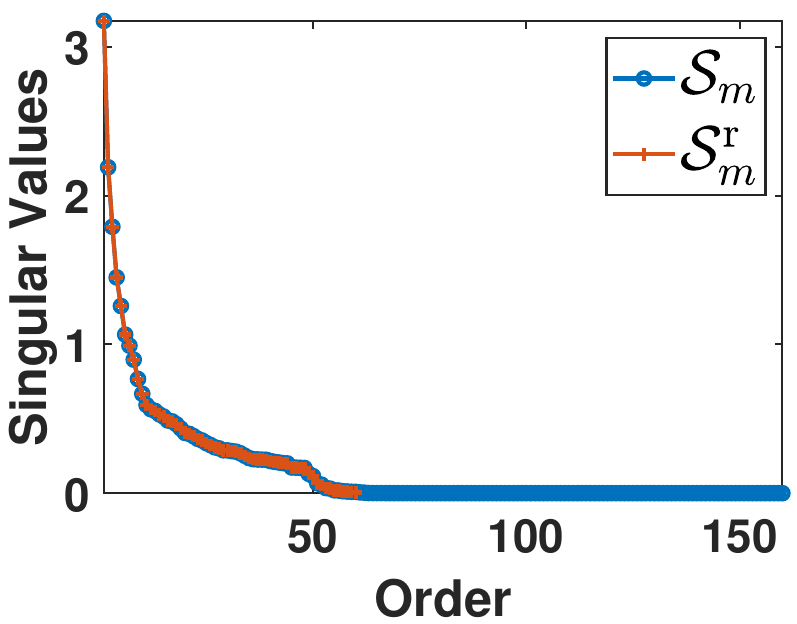}
	\caption{The singular decay of the restricted local solution operator $\Sc_m$ and its low rank approximation $\Sc_m^\mathrm{r}$ for the RTE (left) and the elliptic equation (right). In RTE~\eqref{eqn:rte} we use heterogeneous collision kernel~\eqref{eqn:collision2} $k(x/\eps,v,v') = \sigma^\eps(x) = \frac{1}{81}\frac{1.1+\cos(4\pi x)}{1.1+\sin(2\pi x/\eps)}$, and in elliptic equation we use $a(x,x/\eps) = \frac{2+1.8\sin(\pi x_1/\eps)}{2+1.8\cos(\pi x_2/\eps)} + \frac{2+\sin(\pi x_2/\eps)}{2+1.8\sin(\pi x_1)}$ with $\eps = 2^{-4}$.}
	\label{fig:operator_svd}
\end{figure}

In Figure~\ref{fig:operator_svd} we present the confined solution operator $\Sc_m$ and its low-rank approximation $\Sc_m^\mathrm{r}$. For the radiative transfer equation, the map is a 2880-by-40 matrix, with the size of each patch being 0.2$\times$[-1,1], with $\Delta x = 1/360$ and $\Delta v = 1/40$.
The random sampling procedure reconstructs it with just 6 samples.
For the elliptic equation, the map is  a 1600-by-160 matrix, and the size of each patch is 1$\times$[0,1] with $\Delta x = 1/40$. The random sampling approximates it well with 60 samples.
The compression rates for these examples are thus 6.7 and 2.7, respectively. See details in~\cite{MR4252068}.

\begin{figure}[htb]
	\centering
	\includegraphics[width = 0.23\textwidth]{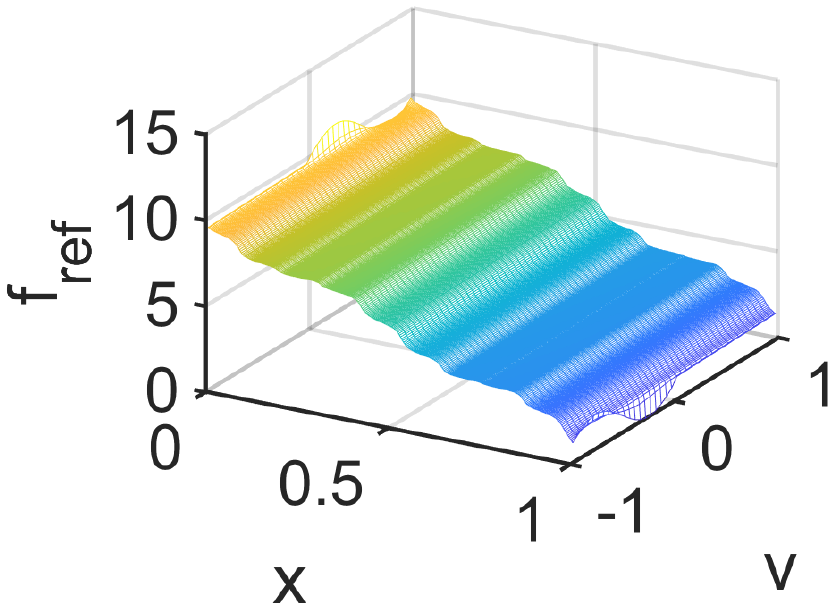}
    \includegraphics[width = 0.23\textwidth]{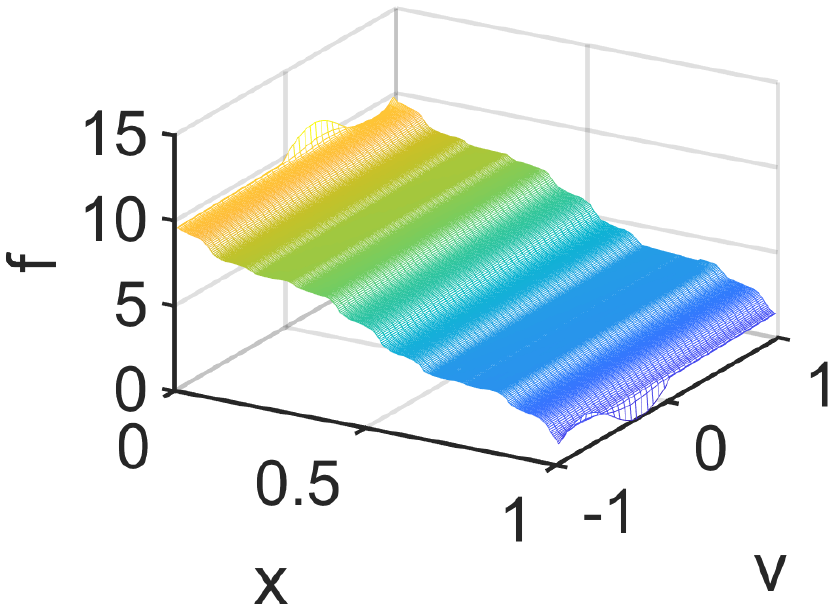}
    \includegraphics[width = 0.23\textwidth]{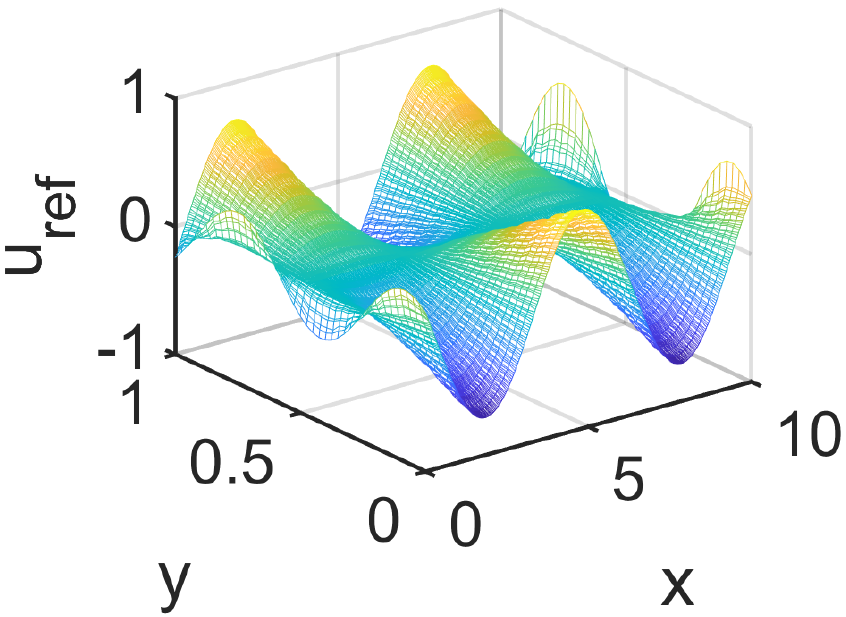}
    \includegraphics[width = 0.23\textwidth]{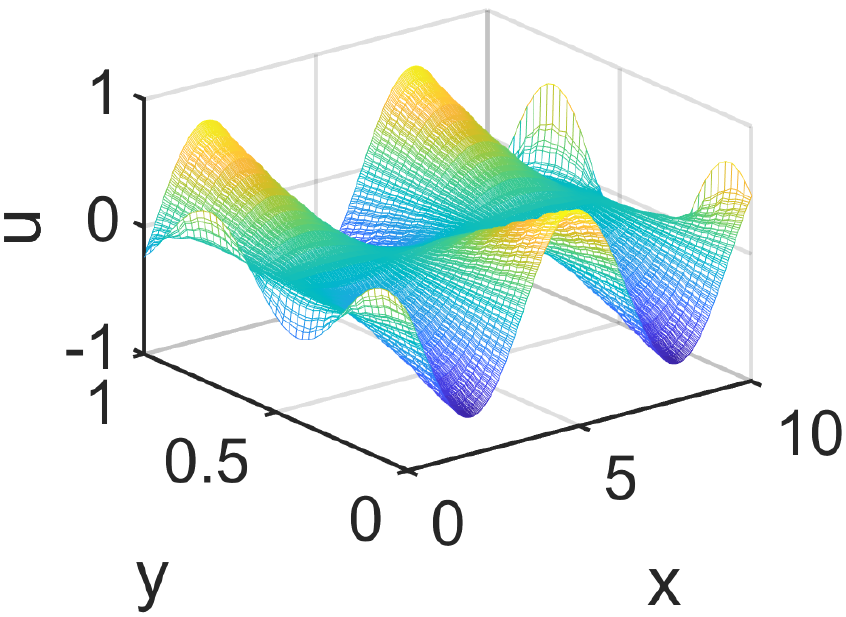}
	\caption{The comparison between the reference solutions (left column) and the approximation using reduced Schwarz method (right column). The first and second rows are for the RTE and the elliptic equation, respectively. In RTE~\eqref{eqn:rte} we use heterogeneous collision kernel~\eqref{eqn:collision2} $k(x/\eps,v,v') = \sigma^\eps(x) = \frac{1}{\eps_1}\frac{1.1+\cos(4\pi x)}{1.1+\sin(2\pi x/\eps_2)}$, with $(\eps_1,\eps_2) = (1/81,1/9)$, and in elliptic equation we use $a(x,x/\eps) = \frac{2+1.8\sin(\pi x_1/\eps)}{2+1.8\cos(\pi x_2/\eps)} + \frac{2+\sin(\pi x_2/\eps)}{2+1.8\sin(\pi x_1)}$ with $\eps = 2^{-4}$.}
	\label{fig:global_Schwarz}
\end{figure}

\begin{figure}[tb]
	\centering
	\includegraphics[width = 0.4\textwidth]{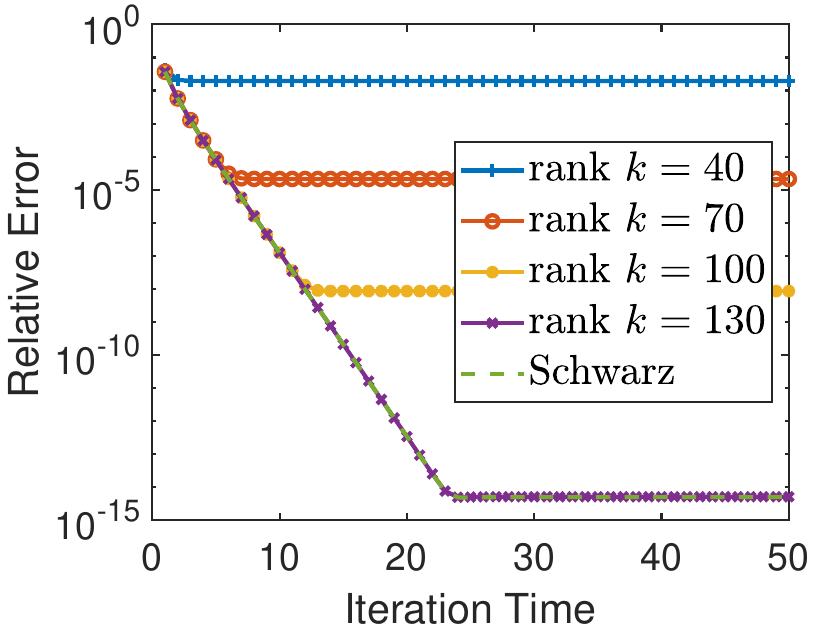}
	\caption{Relative error for reduced Schwarz methods for various ranks $k$ for elliptic equation with oscillatory medium that is the same as Figure~\ref{fig:global_Schwarz}.}
	\label{fig:relative_Schwarz_RTE}
\end{figure}

In Figure~\ref{fig:global_Schwarz}, we show numerical examples for the global solutions of two problems obtained from the approach of this section.
The reference solution (obtained with a fine mesh) is  well captured by the approximation that uses the low-rank BtB map as the surrogate in the Schwarz iteration.
These two cases use just $15\%$ and $43\%$, respectively, of the number of local solves needed to capture the BtB map at fine scale.
While the relative error of the reduced Schwarz method decays as fast as the standard Schwarz iteration, as shown in Figure~\ref{fig:relative_Schwarz_RTE}, the cost is much reduced. See Table~\ref{tbl:time_rte} for a comparison of computation times.

\begin{algorithm}[ht]
	\caption{Reduced Schwarz method for multiscale PDE $\Lc^\eps u^\eps = 0$ over $\Kc$ with $\Bc u^\eps = \phi$ on $\partial\Kc$}\label{alg:schwarz_red}
	\begin{algorithmic}[1]
		\State Given rank $k$, total iterations $T$.
        \State \textbf{Domain Decomposition}
        \Indent
            \State Partition domain according to~\eqref{eqn:partition}.
        \EndIndent
		\State \textbf{Offline Stage:}
        \Indent
		  \State For all $m$, use Algorithm~\ref{alg:rsvd} to find the rank-$k$ approximation to $\Sc_m$, denoted by $\Us_k^m\Sigma_k^m\Vs_k^{m,\ast}$.
        \EndIndent
		\State \textbf{Online}:
        \Indent
		  \State Initiate $f_m^0(x)$ for each $\partial \Kc_m$, and set $t=0$.
		  \State \textbf{While} $|f_n^{t}-f_n^{t-1}|>\text{TOL}$
            \Indent
		  \State Evaluate $u_m^t = \Us_k^m\Sigma_k^m\Vs_k^{m,\ast} f_m^t$ for each $m$.
		  \State Update $f_{n}^{t+1} = u_m^t$ on $\Kc_m\cap\partial\Kc_n$, $n\in\Ic_m$.
            	\State $t\to t+1$.
            \EndIndent
            \State \textbf{End}
		  \State Solve~\eqref{eqn:elliptic_local} for $u_m^t$ using $f_m^t$ for each $m$\,.
		  \State Assemble global solution $u = \sum_{n=1}^{N} \chi_m u^t_m$.
        \EndIndent
		\State \textbf{Return} $u(x)$.
	\end{algorithmic}
\end{algorithm}

\begin{table}
	\centering
	\begin{tabular}{l c | c  }
		\hline \hline
		RTE &offline (s)& online (s) \\
		\hline
		Rank = 5	& 227.99 &   0.0029\\
		Rank = 6	& 268.46 &   0.0162\\
		Full rank   & 706.99 & 0.0148 \\
        \hline
        Schwarz & |	&1027.40\\
		\hline\hline
		elliptic &offline (s)& online (s) \\
		\hline
		Rank = 40	& 49.7 &   0.049\\
		Rank = 70	& 87.3 &   0.061\\
        \hline
        Schwarz & |	&31.4\\
		\hline\hline
	\end{tabular}
    \caption{Run time comparison between vanilla Schwarz method and the reduced Schwarz method. $k=5,6$ for the RTE and $k=40,70$ for the elliptic equation with Dirichlet boundary condition. The configuration of the media is the same as those in Figure~\ref{fig:global_Schwarz}.}
\label{tbl:time_rte}
\end{table}

\section{Manifold learning and nonlinear multiscale problems}\label{sec:manifold}

It is not straightforward to extend the techniques of the previous section to {\em nonlinear} PDEs.
Despite low-rank properties still holding due to the existence of the limiting equation, the argument based on compressing the Green's matrix no longer holds.
The collection of solutions for different source / boundary terms is not longer a linear subspace, but a solution manifold.

We consider the general nonlinear multiscale problem in the following form
\begin{equation}\label{eqn:nonlinear}
\Nc^\eps u^\eps = f\,,
\end{equation}
where $\Nc^\eps$ is a nonlinear differential operator that depends explicitly on the small parameter $\eps$. The term $f$ can be the source term, boundary conditions or initial conditions. Assume further that the equation has an asymptotic limit
\begin{equation}\label{eqn:nonlinear_homo}
\Nc^\ast u^\ast = f
\end{equation}
as $\eps\to0$, that is, $\|u^\eps-u^\ast\| \to 0$ as $\eps\to0$.
The argument for the linear problem is still applicable: The degrees of freedom required by the classical numerical method for solving~\eqref{eqn:nonlinear} grows rapidly as $\epsilon\to0$, while the existence of the homogenized equation~\eqref{eqn:nonlinear_homo} indicates that only $O(1)$ degrees of freedom should be needed to resolve macro-scale features.
From a manifold perspective, the solutions to~\eqref{eqn:nonlinear} vary in a high dimensional space as $f$ changes, but this manifold is approximated to within distance $O(\eps)$ by another manifold whose dimension is $O(1)$.

Suppose a manifold in a high dimensional space is approximately low-dimensional, can we quickly learn it without paying the high dimensional cost?
We turn to manifold learning for answers to this question.
We are particularly interested in adopting the ideas from the local linear embedding and multi-scale SVD approaches that learn the manifold from observed point clouds, and interpolate the local solution manifold using multiple tangent-space patches, see references in~\cite{chen2020manifold}.

We denote the nonlinear solution map of~\eqref{eqn:nonlinear} by $\Sc^\eps:f\in\Xc\to u^\eps\in\Yc$, which maps the source term or initial / boundary conditions $f(x)$ to the solution of the equation.
In the offline stage, we randomly sample a large number of configurations $f_i$ in $\Xc$, and compute the associated solutions $u^\eps_i = \Sc^\eps f_i \in\Yc$ on fine grids.
These solutions form a point cloud in a high dimensional space $\Yc$. The $\{f_i\}$'s are then subdivided into a number of small neighborhoods, and we construct  tangential approximations to the mapping $\Sc^\eps$ on each neighborhood.
In the online stage, given a new configuration $f$, we identify the small neighborhood to which it belongs and find the corresponding solution by performing linear interpolation.
The overall offline-online strategy in summarized in Algorithm~\ref{alg:nonlinear}.
We stress that some modifications are needed to reduce the cost of implementation.
For example, the algorithm should be combined with domain decomposition (for example, Schwarz iteration) to further confine the computation to local domains, to save computational cost.

\begin{figure}[htbp]
  \centering
  \includegraphics[width=0.32\textwidth]{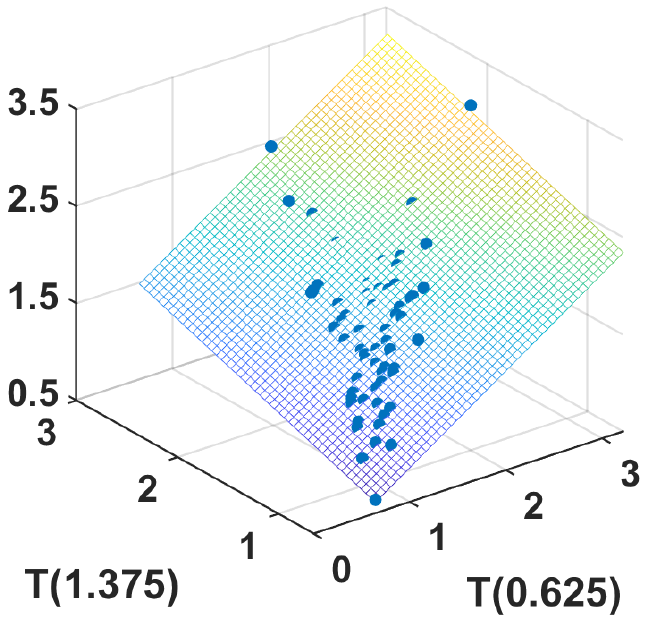}
  \caption{Point cloud and its fitting plane for a 1D nonlinear RTE with Knudsen number $\eps = 2^{-6}$, when confined in a small patch containing an interval $[0.625\,, 1.375]$. The solution profile is approximately determined by the temperature $T$ at two grid points $x=0.625$ and $x=1.375$. The $z$-axis shows the value of $u^\epsilon(x=1)$. The dependence is clearly linear and two-dimensional.}
  \label{fig:offline_manifold}
\end{figure}

In Figure~\ref{fig:offline_manifold} we plot the local low-dimensional solution manifold for a nonlinear RTE (specifically, a linear RTE nonlinearly coupled with a temperature term).
The solution manifold appears to have a local two-dimensional structure; the point clouds lie near on a two-dimensional plane. We refer to~\cite{chen2020manifold} for more details of the implementation and numerical results.

\begin{algorithm}[htbp]
	\caption{Manifold learning algorithm for solving $\Nc^\eps u^\eps = f$}\label{alg:nonlinear}
	\begin{algorithmic}[1]
    \State \textbf{Offline}
	\Indent
		\State Randomly sample $f_i(x)$, $i = 1,\dots,N$, and find solutions $u^\eps_i = \Sc^\eps f_i$.
	\EndIndent
    \State \textbf{Online:} Given $f(x)$:
	\Indent
		\State \textbf{Step 1}: Identify the $k$-nearest neighbors of $f(x)$, call them $f_{i_j}, j = 1,2,\dots,k$, with $f_{i_1}$ being the nearest neighbor;
        \State \textbf{Step 2}: Compute $\Sc^\eps \phi \approx u^\eps_{i_1} +\Us \cdot {c}$ with
    \[
    \Us =
    \begin{bmatrix}
    | & & | \\ u^\eps_{i_2}-u^\eps_{i_1} & \dots & u^\eps_{i_k}-u^\eps_{i_1} \\ | & & |
    \end{bmatrix},
    \]
    where $c$ is a set of coefficient that fits $f-f_{i_1}$ with a
    linear combination of $f_{i_j} - f_{i_1}$, for $j = 2,3,\dotsc,k$.
	\EndIndent
    \State \textbf{Return} $u^\eps = \Sc^\eps \phi$.
    \end{algorithmic}
\end{algorithm}

\section{Looking Forward} \label{sec:conclusion}
We have seen a vast literature addressing all aspects of the computation of multiscale problems. Over the years, the research has been drifting gradually away from its origin, where solvers were influenced by analytical understanding, specifically of the limiting behavior of the specific PDE.
Machine learning algorithms have shown more and more power in sketching the solution profile with a much reduced numerical cost.
In particular, the existence of the homogenized limit suggests there are low-rank features in the discrete system, and that random linear algebra techniques and manifold learning methods, when utilized properly, can identify these features for a compressed representation of the PDE solutions.

We have reviewed two methods, both of which make use of the domain decomposition framework.
They compress either basis functions or the boundary-to-boundary map in an offline learning stage.
This review article serves as a showcase of the power of random solvers in numerical PDEs.
For time-dependent problems, and homogenization problems that have weak-convergence (instead of strong) such as quantum systems in the semi-classical regime, further development of the approaches is needed.
Incorporation of time and weak-limit in the algorithm-design lies at the core of future challenges.

\section{Acknowledgements}
The authors acknowledge generous support from NSF, ONR, DOE, and AFOSR.
Due to restrictions on the allowed number of references, many important contributions are omitted.
The reader is invited to consult the bibliographies in the cited references for a more extensive view of the literature.

\bibliography{ExampleRefs}

\begin{bibdiv}
\begin{biblist}

\bib{MR1185639}{article}{
      author={Allaire, Gr\'{e}goire},
       title={Homogenization and two-scale convergence},
        date={1992},
        ISSN={0036-1410},
     journal={SIAM J. Math. Anal.},
      volume={23},
      number={6},
       pages={1482\ndash 1518},
         url={https://doi-org.ezproxy.library.wisc.edu/10.1137/0523084},
      review={\MR{1185639}},
}

\bib{MR2801210}{article}{
      author={Babuska, Ivo},
      author={Lipton, Robert},
       title={Optimal local approximation spaces for generalized finite element
  methods with application to multiscale problems},
        date={2011},
        ISSN={1540-3459},
     journal={Multiscale Model. Simul.},
      volume={9},
      number={1},
       pages={373\ndash 406},
         url={https://doi-org.ezproxy.library.wisc.edu/10.1137/100791051},
      review={\MR{2801210}},
}

\bib{MR2839402}{book}{
      author={Bensoussan, Alain},
      author={Lions, Jacques-Louis},
      author={Papanicolaou, George~C.},
       title={Asymptotic analysis for periodic structures},
   publisher={AMS Chelsea Publishing, Providence, RI},
        date={2011},
        ISBN={978-0-8218-5324-5},
         url={https://doi-org.ezproxy.library.wisc.edu/10.1090/chel/374},
        note={Corrected reprint of the 1978 original [MR0503330]},
      review={\MR{2839402}},
}

\bib{MR3824169}{article}{
      author={Buhr, Andreas},
      author={Smetana, Kathrin},
       title={Randomized local model order reduction},
        date={2018},
        ISSN={1064-8275},
     journal={SIAM J. Sci. Comput.},
      volume={40},
      number={4},
       pages={A2120\ndash A2151},
         url={https://doi-org.ezproxy.library.wisc.edu/10.1137/17M1138480},
      review={\MR{3824169}},
}

\bib{MR3477310}{article}{
      author={Calo, Victor~M.},
      author={Efendiev, Yalchin},
      author={Galvis, Juan},
      author={Li, Guanglian},
       title={Randomized oversampling for generalized multiscale finite element
  methods},
        date={2016},
        ISSN={1540-3459},
     journal={Multiscale Model. Simul.},
      volume={14},
      number={1},
       pages={482\ndash 501},
         url={https://doi-org.ezproxy.library.wisc.edu/10.1137/140988826},
      review={\MR{3477310}},
}

\bib{MR4155236}{article}{
      author={Chen, Ke},
      author={Li, Qin},
      author={Lu, Jianfeng},
      author={Wright, Stephen~J.},
       title={Random sampling and efficient algorithms for multiscale {PDE}s},
        date={2020},
        ISSN={1064-8275},
     journal={SIAM J. Sci. Comput.},
      volume={42},
      number={5},
       pages={A2974\ndash A3005},
         url={https://doi-org.ezproxy.library.wisc.edu/10.1137/18M1207430},
      review={\MR{4155236}},
}

\bib{chen2020manifold}{article}{
      author={Chen, Shi},
      author={Li, Qin},
      author={Lu, Jianfeng},
      author={Wright, Stephen~J},
       title={Manifold learning and nonlinear homogenization},
        date={2020},
     journal={arXiv preprint arXiv:2011.00568},
}

\bib{MR4252068}{article}{
      author={Chen, Ke},
      author={Li, Qin},
      author={Lu, Jianfeng},
      author={Wright, Stephen~J.},
       title={A {L}ow-{R}ank {S}chwarz {M}ethod for {R}adiative {T}ransfer
  {E}quation {W}ith {H}eterogeneous {S}cattering {C}oefficient},
        date={2021},
        ISSN={1540-3459},
     journal={Multiscale Model. Simul.},
      volume={19},
      number={2},
       pages={775\ndash 801},
         url={https://doi-org.ezproxy.library.wisc.edu/10.1137/19M1276327},
      review={\MR{4252068}},
}

\bib{MR1760042}{article}{
      author={Dumas, Laurent},
      author={Golse, Fran\c{c}ois},
       title={Homogenization of transport equations},
        date={2000},
        ISSN={0036-1399},
     journal={SIAM J. Appl. Math.},
      volume={60},
      number={4},
       pages={1447\ndash 1470},
  url={https://doi-org.ezproxy.library.wisc.edu/10.1137/S0036139997332087},
      review={\MR{1760042}},
}

\bib{MR2830582}{book}{
      author={E, Weinan},
       title={Principles of multiscale modeling},
   publisher={Cambridge University Press, Cambridge},
        date={2011},
        ISBN={978-1-107-09654-7},
      review={\MR{2830582}},
}

\bib{MR1979846}{article}{
      author={E, Weinan},
      author={Engquist, Bjorn},
       title={The heterogeneous multiscale methods},
        date={2003},
        ISSN={1539-6746},
     journal={Commun. Math. Sci.},
      volume={1},
      number={1},
       pages={87\ndash 132},
  url={http://projecteuclid.org.ezproxy.library.wisc.edu/euclid.cms/1118150402},
      review={\MR{1979846}},
}

\bib{MR2477579}{book}{
      author={Efendiev, Yalchin},
      author={Hou, Thomas~Y.},
       title={Multiscale finite element methods},
      series={Surveys and Tutorials in the Applied Mathematical Sciences},
   publisher={Springer, New York},
        date={2009},
      volume={4},
        ISBN={978-0-387-09495-3},
        note={Theory and applications},
      review={\MR{2477579}},
}

\bib{MR1878799}{article}{
      author={Goudon, Thierry},
      author={Mellet, Antoine},
       title={Diffusion approximation in heterogeneous media},
        date={2001},
        ISSN={0921-7134},
     journal={Asymptot. Anal.},
      volume={28},
      number={3-4},
       pages={331\ndash 358},
      review={\MR{1878799}},
}

\bib{MR3445676}{book}{
      author={Hackbusch, Wolfgang},
       title={Hierarchical matrices: algorithms and analysis},
      series={Springer Series in Computational Mathematics},
   publisher={Springer, Heidelberg},
        date={2015},
      volume={49},
        ISBN={978-3-662-47323-8; 978-3-662-47324-5},
  url={https://doi-org.ezproxy.library.wisc.edu/10.1007/978-3-662-47324-5},
      review={\MR{3445676}},
}

\bib{MR3645390}{incollection}{
      author={Hu, Jingwei},
      author={Jin, Shi},
      author={Li, Qin},
       title={Asymptotic-preserving schemes for multiscale hyperbolic and
  kinetic equations},
        date={2017},
   booktitle={Handbook of numerical methods for hyperbolic problems},
      series={Handb. Numer. Anal.},
      volume={18},
   publisher={Elsevier/North-Holland, Amsterdam},
       pages={103\ndash 129},
      review={\MR{3645390}},
}

\bib{MR2806637}{article}{
      author={Halko, Nathan},
      author={Martinsson, Per-Gunnar~J.},
      author={Tropp, Joel~A.},
       title={Finding structure with randomness: probabilistic algorithms for
  constructing approximate matrix decompositions},
        date={2011},
        ISSN={0036-1445},
     journal={SIAM Rev.},
      volume={53},
      number={2},
       pages={217\ndash 288},
         url={https://doi-org.ezproxy.library.wisc.edu/10.1137/090771806},
      review={\MR{2806637}},
}

\bib{Ho:2009multiscale}{article}{
      author={Horstemeyer, Mark~F.},
       title={Multiscale modeling: a review},
        date={2009},
     journal={Practical aspects of computational chemistry},
       pages={87\ndash 135},
}

\bib{MR4050504}{book}{
      author={Martinsson, Per-Gunnar},
       title={Fast direct solvers for elliptic {PDE}s},
      series={CBMS-NSF Regional Conference Series in Applied Mathematics},
   publisher={Society for Industrial and Applied Mathematics (SIAM),
  Philadelphia, PA},
        date={2020},
      volume={96},
        ISBN={978-1-611976-03-8},
      review={\MR{4050504}},
}

\bib{MR3971243}{book}{
      author={Owhadi, Houman},
      author={Scovel, Clint},
       title={Operator-adapted wavelets, fast solvers, and numerical
  homogenization},
      series={Cambridge Monographs on Applied and Computational Mathematics},
   publisher={Cambridge University Press, Cambridge},
        date={2019},
      volume={35},
        ISBN={978-1-108-48436-7},
        note={From a game theoretic approach to numerical approximation and
  algorithm design},
      review={\MR{3971243}},
}

\bib{MR3369060}{article}{
      author={Owhadi, Houman},
       title={Bayesian numerical homogenization},
        date={2015},
        ISSN={1540-3459},
     journal={Multiscale Model. Simul.},
      volume={13},
      number={3},
       pages={812\ndash 828},
         url={https://doi-org.ezproxy.library.wisc.edu/10.1137/140974596},
      review={\MR{3369060}},
}

\bib{MR2382139}{book}{
      author={Pavliotis, Grigorios~A.},
      author={Stuart, Andrew~M.},
       title={Multiscale methods},
      series={Texts in Applied Mathematics},
   publisher={Springer, New York},
        date={2008},
      volume={53},
        ISBN={978-0-387-73828-4},
        note={Averaging and homogenization},
      review={\MR{2382139}},
}

\end{biblist}
\end{bibdiv}

\end{document}